\documentclass[10pt]{article}
\usepackage[latin1]{inputenc}
\usepackage{epsfig}
\usepackage{color}
\usepackage[british,english]{babel}
\usepackage{amsthm}
\usepackage{amsmath}
\usepackage{amsfonts}
\usepackage{amssymb}
\usepackage{graphicx}
\newtheorem{corollary}{Corollary}[section]

\newtheorem{lemma}[corollary]{Lemma}

\newtheorem{remark}[corollary]{Remark}
\newtheorem{theorem}[corollary]{Theorem}
\newfont{\sBlackboard}{msbm10 scaled 900}

\newcommand{\mylabel}[1]{\label{#1}
            \ifx\undefined\stillediting
            \else \fbox{$#1$}\fi }
\newcommand{\BE}{\begin{equation}}

\newcommand{\EEQ}{\end{equation}}
\newcommand{\rfb}[1]{\mbox{\rm
   (\ref{#1})}\ifx\undefined\stillediting\else:\fbox{$#1$}\fi}

\newfont{\Blackboard}{msbm10 scaled 1200}

\newfont{\roma}{cmr10 scaled 1200}

\def\CC{\rm \hbox{C\kern-.56em\raise.4ex
         \hbox{$\scriptscriptstyle |$}\kern+0.5 em }}

\def\n{|\kern -.05cm{|}\kern -.05cm{|}}

\def \noame{\noalign{\medskip}}
\newcommand{\mm}    {{\hbox{\hskip 0.5pt}}}

\newcommand{\bluff} {{\hbox{\raise 15pt \hbox{\mm}}}}

\newcommand{\ep}   {\varepsilon}
\usepackage{fancyhdr}

\makeatletter
\def\section{\@startsection {section}{1}{\z@}{-3.5ex plus -1ex minus
    -.2ex}{2.3ex plus .2ex}{\large\bf}}
\makeatother
\def\be{\begin{equation}}
\def\ee{\end{equation}}

\date{ }
\begin{document}
\thispagestyle{empty}
\title{\Large \bf Analysis of the roughness regimes for micropolar fluids via homogenization}\maketitle
\vspace{-2cm}
\begin{center}
Francisco Javier SU\'AREZ-GRAU\footnote{Departamento de Ecuaciones Diferenciales y An\'alisis Num\'erico. Facultad de Matem\'aticas. Universidad de Sevilla. 41012-Sevilla (Spain) fjsgrau@us.es}
 \end{center}


 \renewcommand{\abstractname} {\bf Abstract}
\begin{abstract} 
We study the asymptotic behavior of micropolar fluid flows in a thin domain of thickness $\eta_\varepsilon$ with a periodic oscillating boundary with wavelength $\varepsilon$. We consider the limit when $\varepsilon$ tends to zero and, depending on the limit of the ratio of $\eta_\varepsilon/\varepsilon$,  we prove the existence of three different regimes. In each regime, we derive a generalized Reynolds equation taking into account the microstructure of the roughness.  
\end{abstract}
\bigskip\noindent

\noindent {\small \bf AMS classification numbers:} 76D08, 76A20, 76A05, 76M50, 35B27, 35Q35.  \\
 
\noindent {\small \bf Keywords:} Homogenization; micropolar fluid flow; Reynolds equation; thin-film fluid.   

\section {Introduction}\label{S1}
We study in this paper the effect of small domain irregularities on thin film flows governed by the linearized 3D micropolar equations. In the case of Newtonian fluids governed by the Stokes or Navier-Stokes equations, this problem has been widely studied since Bayada and Chambat \cite{Bayada1} provided, by means of homogenization techniques, a rigorous derivation of the classical 2D Reynolds equation 
\begin{equation}\label{Reynolds_classiq}
{\rm div}\left(-{h^3\over 12\nu}\nabla p+ b\right)=0\,,
\end{equation}
where  $h$ represents the film thickness, $p$ is the pressure,  $\nu$ is the fluid visco\-sity and $b$ is a vectorial function that usually appears from the exterior forces or from the imposed velocity on a part of the boundary. In this sense,  
va\-rious asymptotic Reynolds-like models, in special regimes, have been obtained depending on the ratio between the size of the roughness and the thickness of the domain and the boundary conditions considered on a part of the boundary, see for example  Bayada {\it et al.} \cite{BCJ}, Benhaboucha {\it et al.} \cite{Benhaboucha}, Benterki {\it et al.} \cite{Benterki}, Bresch {\it et al.} \cite{BCCM}, Boukrouche and Ciuperca \cite{BCiu}, Chupin and Martin \cite{Chupin}, Letoufa {\it et al.} \cite{Letoufa},   Su\'arez-Grau \cite{grau1}, and references therein. 

More precisely, a very general result was obtained in Bayada and Chambat \cite{Bayada_Chambat_2,Bayada_Chambat}, see also Mikelic \cite{Mikelic2}. Assuming that the thickness of the domain is rapidly oscillating, i.e.  the thickness is given by a small parameter $\eta_\ep$ and one of the boundary is rough with small roughness of  wavelength $\ep$, it was proved that depending on the limit of the ratio $\eta_\ep /\ep$, denoted as $\lambda$, there exist three characteristic regimes: Stokes roughness ($0<\lambda<+\infty$), Reynolds roughness ($\lambda=0$) and high-frequency roughness ($\lambda=+\infty$).  In particular, it was obtained that the flow is governed by a generalized 2D Reynolds equation of the form 
\begin{equation}\label{intro_reynolds}{\rm div}\left(-A_\lambda\nabla p+ b_\lambda\right)=0,\end{equation}
for $0\leq \lambda\leq +\infty$, where $A_\lambda$ and $b_\lambda$ are macroscopic quantities known as flow factors, which take into account the microstructure of the roughness. Moreover, it holds that in the Stokes roughness regime the flow factors are calculated by solving 3D local Stokes-like problems depending on the parameter $\lambda$, while in the Reynolds roughness regime they are obtained by solving 2D local Reynolds-like problems, which represents a considerable simplification. In the high-frequency roughness regime,  due to the highly oscillating boundary, the velocity vanishes in the oscillating zone and then, the classical Reynolds equation (\ref{Reynolds_classiq}) is deduced in the non-oscillating zone, so there are no local problems to solve.

 This result has been formally generalized to the unstationary case (the rough surface is moving) in Fabricius {\it et al.} \cite{Fab1}, and recently rigorously gene\-ralized to the case of non-Newtonian fluids governed by the 3D Navier-Stokes system with a nonlinear viscosity (power law)  in Anguiano and Su\'arez-Grau \cite{Anguiano_SG}.

On the other hand, we remark that there are not many papers in the exis\-ting literature dealing with the mathematical modeling of micropolar fluid film lubrication. A generalized version of the Reynolds equation, formally obtained in a critical case when one of the non-Newtonian characteristic parameters has specific (small) order of magnitude, can be found in Singh and Sinha \cite{Sinha} where the authors consider a specific slider-type bearing. Later, in Bayada and Lukaszewicz \cite{BayadaLuc}, it was developed the rigorous derivation, obtaining the generalized version of the 2D  Reynolds equation (\ref{Reynolds_classiq}) for micropolar thin film fluids, which has the form
\begin{equation}\label{Reynolds_classiq_micr}
{\rm div}\left(-{h^3\over 1-N^2}\Phi(h,N)\nabla  p+ b\right)=0,
\end{equation}
where $N$ is the coupling number and 
$$\Phi(h,N)={1\over 12}+{1\over 4h^2(1-N^2)}-{1\over 4h}\sqrt{{N^2\over 1-N^2}}\coth\left(Nh\sqrt{1-N^2}\right).$$

 We also refer to Dupuy {\it et al.} \cite{Dupuy},  for the case of micropolar flow in a curved channel, and  to Marusic-Paloka {\it et al.} \cite{MPM},  for the asymptotic Brinkman-type model proposed starting from 3D micropolar equations.

We remark that in previous papers, the micropolar fluid film  has been considered in a simple thin domain with no roughness introduced. Recently, the roughness effects on a thin film flow have been studied as well and new mathematical models have been proposed in Boukrouche and Paoli \cite{BP}, where the authors consider micropolar flow in a 2D domain assuming the roughness is of the same small order as the film thickness. Employing two-scale convergence technique, they derive the limit problem describing the macroscopic flow. Later, in Pazanin and Su\'arez-Grau \cite{PSG}, a  version of the Reynolds equation is derived in the case of a 3D domain with a particular roughness pattern, where the wavelength of the roughness is assumed to be smaller than the thickness, through a variant of the notion of two-scale convergence introduced in Bresch {\it et al.} \cite{BCCM}.

Our goal in this paper is to give  a general classification result for thin film flows of micropolar fluids with rapidly oscillating thickness in the spirit of  \cite{Anguiano_SG,Bayada_Chambat_2,Bayada_Chambat}, by considering a 3D domain with a thickness given by the parameter $\eta_\ep$ and the wavelength of the roughness by $\ep$. To do this, we use extension results for thin domains and an adaptation of the unfolding method (see Cioranescu {\it et. al} \cite{Ciora,Ciora2}) developed in \cite{Anguiano_SG}. As a result, we deduce that the three characteristic regimes fluids are still valid for micropolar fluids, and moreover,  we derive a generalized version of the Reynolds equation of the form (\ref{intro_reynolds}) depending on $\lambda$. Also, the flow factors are calculated in a different way depending on the regime. More precisely, in the Stokes roughness regime ($0<\lambda<+\infty$) the flow factors are calculated by solving 3D local micropolar Stokes-like problems depending on the parameter $\lambda$, while in the Reynolds roughness regime ($\lambda=0$) they are obtained by solving 2D local micropolar Reynolds-like problems. Finally, in the high-frequency roughness regime ($\lambda=+\infty$) due to the highly oscillating boundary, the classical micropolar Reynolds equation (\ref{Reynolds_classiq_micr}) is deduced in the non-oscillating zone, and there are no local problems to solve.

The paper is organized as follows. In Section \ref{sec:setting} we introduce the domain and some useful notation, and we state the problem. In Section \ref{sec:estimates}, we give some a priori estimates for the velocity,  microrotation and pressure, and we introduce the extension results and the version of the unfolding method necessary to pass to the limit depending on each regime.  The Stokes roughness regime is considered in Section \ref{SRR}, the  Reynolds roughness regime in Section \ref{RRR}, and the high-frequency roughness regime in Section \ref{HFRR}. The corresponding main concergence results are stated in Theorems \ref{them_main_crit}, \ref{them_main_sub} and \ref{thm_general_sup}, respectively.  The paper ends with a conclusion section, an appendix, where we recall the computation of the coefficients of the classical micropolar Reynolds equation (\ref{Reynolds_classiq_micr}),  and with a list of references.

\section{Statement of the problem}\label{sec:setting}
In this section, we first define the thin domain and some sets necessary to study the asymptotic behavior of the solutions.  Next, we introduce the problem considered in the thin domain and also, the rescaled problem posed in a domain of fixed height. We finish this section giving the equivalent weak variational formulation for both problems.

\paragraph{ The domain. }  A  thin domain with a rapidly oscillating thickness is defined by a domain $\omega$ and an associated microstructure given by a function $h_\ep(x')=\eta_\ep h\left(x'/\varepsilon\right)$ that models the roughness of the upper surface and depends  on two small positive parameters $\eta_\ep$ and $\ep$,  representing the thickness of the domain and the wavelength of the roughness, respectively. More precisely, we assume that $\omega$ is an open, smooth, bounded and connected set of $\mathbb{R}^2$, and $h$ is a positive and smooth function, defined for $y'$ in $\mathbb{R}^2$, $Y'$-periodic with $Y'=(-1/2,1/2)^2$ the cell of periodicity in $\mathbb{R}^2$, and there exist $h_{\rm min}$ and $h_{\rm max}$ such that
$$0<h_{\rm min}=\min_{y'\in Y'} h(y'),\quad h_{\rm max}=\max_{y'\in Y'}h(y')\,.$$
We remark that along this paper, the points $x\in \mathbb{R}^3$ will be decomposed as $x=(x',x_3)$ with $x'\in\mathbb{R}^2$, $x_3\in\mathbb{R}$. We also use the notation $x'$ to denote a generic vector of $\mathbb{R}^2$.

Thus, we define the thin domain $\Omega_\ep\subset\mathbb{R}^3$ by 
$$\Omega_\ep=\left\{(x',x_3)\in \mathbb{R}^2\times \mathbb{R}\,:\,x'\in\omega,\ 0<x_3< h_\ep(x')\right\}\,,$$
and the oscillating part of the boundary  by $\Sigma_\varepsilon=\omega\times \{h_\ep(x')\}$. We also define the respective rescaled sets $\widetilde\Omega_\ep=\omega\times (0,h(x'/\ep))$ and  $\widetilde \Sigma_\varepsilon=\omega\times\{h(x'/\ep)\}$.

Related  to the microstructure of the periodicity of the  boundary,  we consider that the domain $\omega$ is covered by a rectangular mesh of size $\ep$: for $k'\in\mathbb{Z}^2$, each cell $Y'_{k',\ep}=\ep k'+\ep Y'$, and for simplicity, we assume that there exists an exact finite number of periodic sets $Y'_{k',\ep}$.  We define $T_\ep=\{k'\in\mathbb{Z}^2\,:\, Y'_{k',\ep}\cap\omega\neq \emptyset\}$. Also, we define $Y_{k',\ep}=Y'_{k',\ep}\times (0,h(y'))$ and  $Y=Y'\times (0,h(y'))$, which is the reference cell in $\mathbb{R}^3$.

Two quantities  $h_{\rm min}$ and $h_{\rm max}$ allow us to define:
\begin{itemize}
\item[--] The extended sets $Q_\ep=\omega\times (0,\eta_\varepsilon h_{\rm max})$,  $\Omega=\omega\times  (0, h_{\rm max})$  and $\Sigma=\omega\times \{h_{\rm max}\}$.
\item[--]  The extended cube  $\widetilde Q_{k',\varepsilon}=Y'_{k',\varepsilon}\times (0, h_{\rm max})$ for $k'\in\mathbb{Z}^2$.

\item[--] The restricted sets $\Omega_\ep^+=\omega\times (\eta_\varepsilon h_{\rm min},h_\ep(x'))$, $\widetilde \Omega_\ep^+=\omega\times (h_{\rm min},h(x'/\ep))$, $\Omega^+=\omega\times (h_{\rm min}, h_{\rm max})$ and $\Omega^-=\omega\times (0,h_{\rm min})$.

\item[--]  The extended and restricted basic cells $\Pi=Y'\times (0,h_{\rm max})$,  $\Pi^+=Y'\times (h_{\rm min},h_{\rm max})$ and $\Pi^-=Y'\times (0,h_{\rm min})$.
\end{itemize}
In order to apply the unfolding method, we will use the following notation. For $k'\in \mathbb{Z}^2$, we define $\kappa: \mathbb{R}^2\to \mathbb{Z}^2$ by
\begin{equation}\label{kappa_fun}
\kappa(x')=k' \Longleftrightarrow x'\in Y'_{k',1}\,.
\end{equation}
Remark that $\kappa$ is well defined up to a set of zero measure in $\mathbb{R}^2$ (the set $\cup_{k'\in \mathbb{Z}^2}\partial Y'_{k',1}$). Moreover, for every $\varepsilon>0$, we have
$$\kappa\left({x'\over \varepsilon}\right)=k'\Longleftrightarrow x'\in Y'_{k',\varepsilon}\,.$$
We denote by $O_\ep$ a generic real sequence which tends to zero with $\ep$ and can change from line to line. We denote by $C$ a generic constant which can change from line to line. To finish, let $C^\infty_{\#}(Y)$ be the space of infinitely differentiable functions in $\mathbb{R}^3$ that are $Y'$-periodic. By $L^2_{\#}(Y)$ (resp. $H^1_{\#}(Y)$) we denote its completion in the norm $L^2(Y)$ (resp. $H^1(Y)$) and by $L^2_{0,\#}(Y)$  the space of functions in $L^2_{\#}(Y)$ with zero mean value.

\paragraph{The problem.} When the distance between two surfaces becomes very small, the experimental results from the tribology literature (see e.g. \cite{John,Luo,Luo2}) suggest that the fluid's internal structure should be taken into account as well. 
Among various non-Newtonian models, the model of micropolar fluid (proposed by Eringen \cite{Eringen}
in 60's) turns out to be the most appropriate since it acknowledges
the effects of the local structure and micro-motions of the fluid
elements. Physically, micropolar fluids consist in a large number of
small spherical particles uniformly dispersed in a viscous medium.
Assuming that the particles are rigid and ignoring their deformations,
the related mathematical model expresses the balance of
momentum, mass and angular momentum. A new unknown function
called microrotation (i.e. the angular velocity field of rotation
of particles) is added to the usual velocity and pressure fields. Consequently,
Navier-Stokes equations become coupled with a new
vector equation coming from the conservation of angular momentum
with four microrotation viscosities introduced (see \cite{Luka} for more details). Being able to
describe nume\-rous real fluids better than the classical (Newtonian)
model, micropolar fluid models have been extensively studied in
recent years (see e.g. \cite{Bonnivard,BP,Dupuy,PSG}).

Taking into account the application we want to model (lubrication with micropolar fluid), it
is reasonable to assume a small Reynolds number and omit the inertial terms in momentum equations of the micropolar system. Also, it has been observed  that the magnitude
of the viscosity coefficients appearing in the micropolar equations may influence
the effective flow. Thus, it is reasonable to work with the system
written in a non-dimensional form (see e.g. \cite{BayadaLuc} for more details). Thus,  we consider the stationary flow of an incompressible micropolar fluid in $\Omega_\ep$ which is governed  by the following linearized micropolar system formulated in a non-dimensional form
\begin{equation}\label{system_1}
\left\{\begin{array}{rl}
-{\rm div}(D u_\ep)+\nabla p_\ep=2N^2{\rm rot}\,w_\ep+ f_\ep&\quad\hbox{in}\quad\Omega_\ep,\\
\noame
{\rm div}\,u_\ep=0&\quad\hbox{in}\quad\Omega_\ep,\\
\noame
-R_M{\rm div}(D w_\ep)+4N^2w_\ep=2N^2{\rm rot}\,u_\ep+g_\ep&\quad\hbox{in}\quad\Omega_\ep,
\end{array}\right.
\end{equation}
with homogeneous boundary conditions (it does not alter the generality of the problem under consideration),
\begin{equation}\label{bc_system_1}
u_\ep=w_\ep=0\quad \hbox{on}\quad\partial\Omega_\ep\,.
\end{equation}
In system (\ref{system_1})-(\ref{bc_system_1}), the velocity $u_\ep$, the pressure $p_\ep$ and the microrotation $w_\ep$ are unknown. Dimensionless (non-Newtonian) parameter $N^2$ characterizes the coupling between the equations for the velocity and microrotation and it is of order $\mathcal{O}(1)$, in fact $N^2 $ lies between zero and one. The second dimensionless parameter, denoted by $R_M$ is, in fact, related to the characteristic length of the microrotation effects and is compared with the small parameter $\eta_\ep$ by assuming that $R_M=\mathcal{O}(\eta_\ep^2)$, namely
\begin{equation}\label{R_M}
R_M=\eta_\ep^2R_c\quad\hbox{with }R_c=\mathcal{O}(1)\,.
\end{equation}
This case is the situation that is commonly introduced to study the micropolar fluid because the third equation of (\ref{system_1}) shows a strong coupling between velocity and microrotation in the limit (see   \cite{BayadaChamGam,BayadaLuc}).

Under assumptions that $f_\ep,g_\ep\in L^2(\Omega_\ep)^3$, it is well known that problem (\ref{system_1})-(\ref{bc_system_1}) has a unique weak solution $(u_\ep,w_\ep,p_\ep)\in H^1_0(\Omega_\ep)^3\times H^1_0(\Omega_\ep)^3\times L^2_0(\Omega_\ep)$ (see  \cite{Luka}), where the space $L^2_0$ is the space of functions of $L^2$ with null integral.

Our aim is to study the asymptotic behavior of $u_\ep$, $w_\ep$ and $p_\ep$ when $\ep$ and $\eta_\ep$ tend to zero and identify homogenized models coupling the effects of the thickness of the domain and the roughness of the boundary. For this purpose, as usual when we deal with thin domains, we use the dilatation in the variable $x_3$ given by
\begin{equation}\label{dilatacion}
y_3={x_3\over \eta_\ep}\,,
\end{equation}
in order to have the functions defined in the open set with fixed height   $\widetilde\Omega_\ep$ with oscillating boundary $\widetilde \Sigma_\varepsilon$.

Namely, we define  $\tilde u_\ep, \tilde w_\ep\in H^1_0(\widetilde\Omega_\ep)^3$ and $\tilde p_\ep\in L^2_0(\widetilde\Omega_\ep)$ by 
\begin{equation}\label{unk_dilat}
\begin{array}{c}
\tilde u_\ep(x',y_3)=u_\ep(x',\eta_\ep y_3),\quad \tilde w_\ep(x',y_3)=w_\ep(x',\eta_\ep y_3),\\
\noame \tilde p_\ep(x',y_3)=p_\ep(x',\eta_\ep y_3),\quad\hbox{a.e. }(x',y_3)\in \widetilde\Omega_\ep\,.
\end{array}
\end{equation}
Let us introduce some notation which will be useful in the following. For a vectorial function $v=(v',v_3)$ and a scalar function $w$, we introduce the operators $D_{\eta_\ep}$, $\nabla_{\eta_\ep}$ and ${\rm rot}_{\eta_\ep}$ by 
\begin{eqnarray}
&(D_{\eta_\ep}v)_{ij}=\partial_{x_j}v_i\hbox{ for }i=1,2,3,\ j=1,2,\quad (D_{\eta_\ep})_{i,3}={1\over \eta_\ep}\partial_{y_3}v_i\hbox{ for }i=1,2,3,\nonumber&\\
\noame
&\nabla_{\eta_\ep}w=(\nabla_{x'}w,{1\over \eta_\ep}\partial_{y_3}w)^t,\quad {\rm div}_{\eta_\ep}v={\rm div}_{x'}v'+{1\over \eta_\ep}\partial_{y_3}v_3,\nonumber&\\
\noame
&{\rm rot}_{\eta_\ep}v=\left({\rm rot}_{x'}v_3+{1\over \eta_\ep}{\rm rot}_{y_3}v',{\rm Rot}_{x'}v'\right)^t,\nonumber&
\end{eqnarray}
where, denoting $(v')^\perp=(-v_2,v_1)^t$, we define
\begin{equation}\label{def_rot}
\begin{array}{c}
{\rm rot}_{x'}v_3=(\partial_{x_2}v_3,-\partial_{x_1}v_3)^t,\quad {\rm rot}_{y_3}v'=\partial_{y_3}(v')^\perp,\\
\noame {\rm Rot}_{x'}v'=\partial_{x_1}v_2-\partial_{x_2}v_1.
\end{array}
\end{equation}
Using the transformation (\ref{dilatacion}), the rescaled system (\ref{system_1})-(\ref{bc_system_1}) can be rewritten as
\begin{equation}\label{system_2}
\left\{
\begin{array}{rl}
-{\rm div}_{\eta_\ep}(D_{\eta_\ep} \tilde u_\ep)+\nabla_{\eta_\ep} \tilde p_\ep=2N^2{\rm rot}_{\eta_\ep}\,\tilde w_\ep+ \tilde f_\ep&\quad\hbox{in}\quad\widetilde\Omega_\ep,\\
\noame
{\rm div}_{\eta_\ep}\tilde u_\ep=0&\quad\hbox{in}\quad\widetilde \Omega_\ep,\\
\noame
-\eta_\ep^2R_c{\rm div}_{\eta_\ep}(D_{\eta_\ep} \tilde w_\ep)+4N^2\tilde w_\ep=2N^2{\rm rot}_{\eta_\ep}\tilde u_\ep+\tilde g_\ep&\quad\hbox{in}\quad\widetilde\Omega_\ep\,,
\end{array}\right.
\end{equation}
with homogeneous boundary conditions
\begin{equation}\label{bc_system_2}
\tilde u_\ep=\tilde w_\ep=0\quad \hbox{on}\quad\partial\widetilde\Omega_\ep\,,
\end{equation}
where $\tilde f_\ep$ and $\tilde g_\ep$ are defined similarly as in (\ref{unk_dilat}). 

Our goal then is to describe the asymptotic behavior of this new sequences $\tilde u_\ep$, $\tilde w_\ep$ and $\tilde p_\ep$ when $\ep$ and $\eta_\ep$ tend to zero. To do this, we need to obtain appropriate estimates, so it is usual to consider for $f_\ep$  and $g_\ep$  the following estimates 
\begin{equation}\label{estim_f_g_ep}
\|f_\ep\|_{L^2(\Omega_\ep)^3}\leq C\eta_\ep^{1\over 2},\quad \|g_\ep\|_{L^2(\Omega_\ep)^3}\leq C\eta_\ep^{3\over 2}\,.
\end{equation}
For example, assuming $f,g\in L^2(\Omega)$, we can consider as external forces satisfying (\ref{estim_f_g_ep})  the following ones
$$f_\ep(x)=f\left(x',{x_3\over \eta_\ep}\right),\quad g_\ep(x)=\eta_\ep g\left(x',{x_3\over \eta_\ep}\right)\quad\hbox{a.e. }x\in\Omega_\ep.$$
We point out that due to the thickness of the domain, it is usual to assume that the vertical components of the external forces can be neglected and, moreover the forces can be considered independent of the vertical variable. Thus, for sake of simplicity, assuming $f',g'\in L^2(\omega)^2$, along the paper we will consider the following assumptions:

\begin{itemize}
\item[(i)] If $\eta_\ep\approx\ep$, with $\eta_\ep/\ep\to \lambda$, $0<\lambda<+\infty$, or $\eta_\ep\ll\ep$, then
\begin{equation}\label{estim_f_g_cases_crit_sub}
f_\ep(x)=(f'(x'),0),\quad g_\ep=(\eta_\ep g'(x'),0),\quad\hbox{a.e. }x\in\Omega_\ep.
\end{equation}
\item[(ii)] If $\eta_\ep\gg\ep$, then
\begin{equation}\label{estim_f_g_cases_sup}
f_\ep(x)=(f'(x'),0),\quad g_\ep=(\ep g'(x'),0),\quad\hbox{a.e. }x\in\Omega_\ep.
\end{equation}
\end{itemize}
We observe that in this case $\tilde f_\ep=f_\ep$ and $\tilde g_\ep=g_\ep$ and that in $(i)$ the external forces satisfy (\ref{estim_f_g_ep}).  However, in the case $(ii)$, due to the high oscillations of the boundary, to obtain appropriate estimates it is necessary to assume that $g_\ep$ satisfies a more precise estimate, that is
$\| g_\ep\|_{L^2(\Omega_\ep)^3}\leq C \ep\eta_\ep^{1\over 2}$ (see proof of Lemma \ref{lemma_estimates} for more details). 

\paragraph{Weak variational formulations.} We finish this section by giving the equivalent weak variational formulation of  system (\ref{system_1})-(\ref{bc_system_1}) and  the rescaled system (\ref{system_2})-(\ref{bc_system_2}), which will be useful in next sections.  

For  problem (\ref{system_1})-(\ref{bc_system_1}), the weak variational formulation is to find $u_\ep, w_\ep\in H^1_0(\Omega_\ep)^3$ and $p_\ep\in L^2_0(\Omega_\ep)$ such that
\begin{equation}\label{form_var_1}
\left\{\begin{array}{l}
\displaystyle \int_{\Omega_\ep}D u_\ep:D\varphi\,dx-\int_{\Omega_\ep}p_\ep\,{\rm div}\,\varphi\,dx\\
\noame
\displaystyle=2N^2\int_{\Omega_\ep}{\rm rot}\,w_\ep\cdot \varphi\,dx+\int_{\Omega_\ep}f_\ep\cdot\varphi\,dx,\\
\\
\displaystyle \eta_\ep^2R_c \int_{\Omega_\ep}D w_\ep:D\psi\,dx+4N^2\int_{\Omega_\ep}w_\ep\cdot\psi\,dx\\
\noame
\displaystyle=2N^2\int_{\Omega_\ep}{\rm rot}\,u_\ep\cdot \psi\,dx+\int_{\Omega_\ep}g_\ep\cdot\psi\,dx\,,
\end{array}\right.
\end{equation}
for every $\varphi,\psi \in H^1_0(\Omega_\ep)^3$, and  the equivalent  weak variational formulation for the rescaled system (\ref{system_2})-(\ref{bc_system_2}) is to find $\tilde u_\ep, \tilde w_\ep\in H^1_0(\widetilde \Omega_\ep)^3$ and $\tilde p_\ep\in L^2_0(\widetilde \Omega_\ep)$ such that
\begin{equation}\label{form_var_2}
\left\{\begin{array}{l}
\displaystyle \int_{\widetilde \Omega_\ep}D_{\eta_\ep} \tilde u_\ep:D_{\eta_\ep}\tilde \varphi\,dx'dy_3-\int_{\widetilde \Omega_\ep}\tilde p_\ep\,{\rm div}_{\eta_\ep}\tilde \varphi\,dx'dy_3\\
\noame
\displaystyle=2N^2\int_{\widetilde \Omega_\ep}{\rm rot}_{\eta_\ep}\tilde w_\ep\cdot\tilde  \varphi\,dx'dy_3+\int_{\widetilde \Omega_\ep}\tilde f_\ep\cdot\tilde \varphi\,dx'dy_3\,,\\
\\
\displaystyle \eta_\ep^2R_c \!\!\int_{\widetilde \Omega_\ep}\!D_{\eta_\ep} \tilde w_\ep:D_{\eta_\ep}\tilde \psi\,dx'dy_3+4N^2\!\int_{\widetilde \Omega_\ep}\!\!\tilde w_\ep\cdot\tilde \psi\,dx'dy_3\\
\noame
\displaystyle=2N^2\!\int_{\widetilde \Omega_\ep}\!{\rm rot}_{\eta_\ep}\tilde u_\ep\cdot \tilde \psi\,dx'dy_3+\!\int_{\widetilde \Omega_\ep}\!\!\tilde g_\ep\cdot\tilde \psi\,dx'dy_3\,,
\end{array}\right.
\end{equation}
for every $\tilde \varphi,\tilde \psi \in H^1_0(\widetilde \Omega_\ep)^3$.

\section{{\it A priori} estimates}\label{sec:estimates}
First, we recall the Poincar\'e inequality in a domain with thickness $\eta_\ep$ (see \cite{MT}).
\begin{lemma}
For every $v\in H^1_0(\Omega_\ep)^3$, the following inequality holds
\begin{equation}\label{Poincare}
\|v\|_{L^2(\Omega_\ep)^3}\leq c_2\eta_\ep\|Dv\|_{L^2(\Omega_\ep)^{3\times 3}}\,,
\end{equation}
where $c_2>0$ is independent of $v$, $\ep$ and $\eta_\ep$.
\end{lemma}
Next, we give the following results relating the derivative and the rotational.
\begin{lemma}
For $v\in H^1_0(\Omega_\ep)^3$, the following inequality holds
\begin{equation}\label{Gaffney}
\|{\rm rot}\,v\|_{L^2(\Omega_\ep)^{3}}\leq \|Dv\|_{L^2(\Omega_\ep)^{3\times 3}},
\end{equation}
and, if moreover ${\rm div}\,v=0$ in $\Omega_\ep$, then it holds
\begin{equation}\label{Gaffney_div0}
\|{\rm rot}\,v\|_{L^2(\Omega_\ep)^{3}}=\|Dv\|_{L^2(\Omega_\ep)^{3\times 3}}.
\end{equation}
\end{lemma}
\begin{proof} By using relation $-\Delta v ={\rm rot}\,({\rm rot}\, v)-\nabla\,{\rm div}\varphi$,  it can be proved (see \cite{DuvautLions}) that
$$\int_{\Omega_\ep}|Dv|^2dx=\int_{\Omega_\ep}|{\rm rot}\,v|^2dx+ \int_{\Omega_\ep}|{\rm div}\,v|^2dx,\quad \forall\,v\in H^1_0(\Omega_\ep)^3.$$
Then, (\ref{Gaffney}) easily holds, and (\ref{Gaffney_div0}) is a consequence of the free divergence condition.
\hfill$\square$
\end{proof}
We start by obtaining some {\it a priori} estimates for $u_\ep$, $w_\ep$, $\tilde u_\ep$ and $\tilde w_\ep$.

\begin{lemma}\label{lemma_estimates}
There exists a constant $C$ independent of $\ep$, such that the solution $(u_\ep,w_\ep)$ of problem (\ref{system_1})-(\ref{bc_system_1}) and the corresponding rescaled solution $(\tilde u_\ep,\tilde w_\ep)$ of the problem (\ref{system_2})-(\ref{bc_system_2}) satisfy
\begin{eqnarray}
&\|u_\ep\|_{L^2(\Omega_\ep)^3}\leq C\eta_\ep^{5\over 2},&\quad \|Du_\ep\|_{L^2(\Omega_\ep)^{3\times 3}}\leq C\eta_\ep^{3\over 2}\,,\label{estim_u_ep}\\
\noame
&\|w_\ep\|_{L^2(\Omega_\ep)^3}\leq C\eta_\ep^{3\over 2},&\quad \|Dw_\ep\|_{L^2(\Omega_\ep)^{3\times 3}}\leq C\eta_\ep^{1\over 2}\,,\label{estim_w_ep}
\end{eqnarray}
\begin{eqnarray}
&\|\tilde u_\ep\|_{L^2(\widetilde \Omega_\ep)^3}\leq C\eta_\ep^{2},&\quad \|D_{\eta_\ep}\tilde u_\ep\|_{L^2(\widetilde \Omega_\ep)^{3\times 3}}\leq C\eta_\ep\,,\label{estim_tilde_u_ep}\\
\noame
&\|\tilde w_\ep\|_{L^2(\widetilde \Omega_\ep)^3}\leq C\eta_\ep,&\quad \|D_{\eta_\ep}\tilde w_\ep\|_{L^2(\widetilde \Omega_\ep)^{3\times 3}}\leq C\,.\label{estim_tilde_w_ep}
\end{eqnarray}
Moreover, in the case $\eta_\ep\gg \ep$, defining the restriction functions $u_\ep^+:=u_\ep|_{\Omega_\ep^+}$, $w_\ep^+:=w_\ep|_{\Omega_\ep^+}$, $\tilde u_\ep^+:=\tilde u_\ep|_{\widetilde \Omega_\ep^+}$ and $\tilde w_\ep^+:=\tilde w_\ep|_{\widetilde \Omega_\ep^+}$, we also have the following estimates
\begin{eqnarray}
&\|u_\ep^+\|_{L^2(\Omega_\ep^+)^3}\leq \eta_\ep^{1\over 2}\ep^2,&\quad \|Du_\ep^+\|_{L^2(\Omega_\ep^+)^{3\times 3}}\leq \eta_\ep^{1\over 2}\ep,\label{estim_u_ep+}\\
\noame
&\|w_\ep^+\|_{L^2(\Omega_\ep^+)^3}\leq \eta_\ep^{-{3\over 2}}\ep^3,&\quad \|Dw_\ep^+\|_{L^2(\Omega_\ep^+)^{3\times 3}}\leq \eta_\ep^{-{3\over 2}}\ep^2,\label{estim_w_ep+}
\\
\noame
&\|\tilde u_\ep^+\|_{L^2(\widetilde \Omega_\ep^+)^3}\leq \ep^2,&\quad \|D_{\eta_\ep}\tilde u_\ep^+\|_{L^2(\widetilde \Omega_\ep^+)^{3\times 3}}\leq \ep,\label{estim_tilde_u_ep+}\\
\noame
&\|\tilde w_\ep^+\|_{L^2(\widetilde \Omega_\ep^+)^3}\leq \eta_\ep^{-2}\ep^3,&\quad \|D_{\eta_\ep}\tilde w_\ep^+\|_{L^2(\widetilde \Omega_\ep^+)^{3\times 3}}\leq \eta_\ep^{-2}\ep^2.\label{estim_tilde_w_ep+}
\end{eqnarray}
\end{lemma}
\begin{proof} For every cases, taking $\varphi=u_\ep$ as test function in the first equation of (\ref{form_var_1}), taking into account $\int_{\Omega_\ep}{\rm rot}\,w_\ep\cdot u_\ep\,dx=\int_{\Omega_\ep}{\rm rot}\,u_\ep\cdot w_\ep\,dx$, applying Cauchy-Schwarz's inequality and from (\ref{estim_f_g_cases_crit_sub}), (\ref{estim_f_g_cases_sup}), (\ref{Poincare}) and (\ref{Gaffney_div0}), we have 
\begin{equation}\label{lem_estim_1}
\begin{array}{l}
\displaystyle
\|Du_\ep\|^2_{L^2(\Omega_\ep)^{3\times 3}} \displaystyle =2N^2\int_{\Omega_\ep}{\rm rot}\,w_\ep\cdot u_\ep\,dx+\int_{\Omega_\ep}f_\ep\cdot u_\ep\,dx \\
\noame
\displaystyle 
\quad=2N^2\int_{\Omega_\ep}w_\ep\cdot {\rm rot}\,u_\ep\,dx+\int_{\Omega_\ep}f'(x')\cdot u_\ep'\,dx\\
\noame
 \displaystyle
 \quad \leq 2N^2 \|w_\ep\|_{L^2(\Omega_\ep)^3}\|D u_\ep\|_{L^2(\Omega_\ep)^{3\times 3}}+\eta_\ep^{3\over 2}c_2\|f'\|_{L^2(\omega)^2}\|D u_\ep\|_{L^2(\Omega_\ep)^{3\times 3}}\,, 
\end{array}
\end{equation}
which implies 
\begin{equation}\label{lem_estim_2}
\eta_\ep^{-{3\over 2}}\|D u_\ep\|_{L^2(\Omega_\ep)^{3\times 3}}\leq \eta_\ep^{-{3\over 2}}2N^2\|w_\ep\|_{L^2(\Omega_\ep)^3}+c_2\|f'\|_{L^2(\omega)^2}\,.
\end{equation}
In the cases $\eta_\ep\approx\ep$ and $\eta_\ep\ll\ep$, taking $\psi=w_\ep$ as test function in the second equation of (\ref{form_var_1}), applying Cauchy-Schwarz's inequality and taking into account (\ref{estim_f_g_cases_crit_sub}), we have
\begin{equation}\label{lem_estim_3}
\begin{array}{l}
\displaystyle
\eta_\ep^2R_c\|D w_\ep\|_{L^2(\Omega_\ep)^{3\times 3}}^2+4N^2\|w_\ep\|_{L^2(\Omega_\ep)^3}^2\\
\noame
\displaystyle\quad\quad  =2N^2\int_{\Omega_\ep}{\rm rot}\,u_\ep\cdot w_\ep\,dx+\eta_\ep\int_{\Omega_\ep}g'(x')\cdot w_\ep'\,dx\\
\noame
\displaystyle\quad\quad 
\leq 2N^2\|w_\ep\|_{L^2(\Omega_\ep)^3}\|D u_\ep\|_{L^2(\Omega_\ep)^{3\times 3}}+\eta_\ep^{3\over 2}\|g'\|_{L^2(\omega)^2}\|w_\ep\|_{L^2(\Omega_\ep)^3}\,,
\end{array}
\end{equation}
which implies
\begin{equation}\label{lem_estim_4}
\eta_\ep^{-{3\over 2}}2N^2\|w_\ep\|_{L^2(\Omega_\ep)^3}\leq \eta_\ep^{-{3\over 2}}N^2\|D u_\ep\|_{L^2(\Omega_\ep)^{3\times 3}}+{1\over 2}\|g'\|_{L^2(\omega)^2}\,.
\end{equation}
In the case $\eta_\ep\gg \ep$, proceeding as above by taking into account (\ref{estim_f_g_cases_sup}), and using that in this case
$$\ep\eta_\ep^{1\over 2}\|g'\|_{L^2(\omega)^2}\|w_\ep\|_{L^2(\Omega_\ep)^3}\leq \eta_\ep^{3\over 2}\|g'\|_{L^2(\omega)^2}\|w_\ep\|_{L^2(\Omega_\ep)^3},$$
then estimate (\ref{lem_estim_4}) also holds. 

Then, from (\ref{lem_estim_2}) and (\ref{lem_estim_4}), we conclude for every cases that
$$\eta_\ep^{-{3\over 2}}\|D u_\ep\|_{L^2(\Omega_\ep)^{3\times 3}}\leq {c_2\over 1-N^2}\|f'\|_{L^2(\omega)^2}+{1\over 2(1-N^2)}\|g'\|_{L^2(\omega)^2}\,,$$
which gives the second estimate in (\ref{estim_u_ep}). This together with (\ref{Poincare}) gives the first one. Moreover, by means of the dilatation (\ref{dilatacion}), we get (\ref{estim_tilde_u_ep}).\\

\noindent To get the second estimate in (\ref{estim_w_ep}), we use $\int_{\Omega_\ep}{\rm rot}\,u_\ep\cdot w_\ep\,dx=\int_{\Omega_\ep}{\rm rot}\,w_\ep\cdot u_\ep\,dx$  in (\ref{lem_estim_3}), (\ref{Poincare}) and (\ref{Gaffney}), and proceeding as above we obtain in every cases
\begin{equation}\label{lem_estim_5}
\begin{array}{l}
\displaystyle
\eta_\ep^2R_c\|D w_\ep\|_{L^2(\Omega_\ep)^{3\times 3}}^2+4N^2\|w_\ep\|_{L^2(\Omega_\ep)^3}^2\\
\noame
\displaystyle\quad
\leq 2N^2\|u_\ep\|_{L^2(\Omega_\ep)^3}\|D w_\ep\|_{L^2(\Omega_\ep)^{3\times 3}}+\eta_\ep^{5\over 2}\|g'\|_{L^2(\omega)^2}\|Dw_\ep\|_{L^2(\Omega_\ep)^{3\times 3}}\,,
\end{array}
\end{equation}
which, by using the estimate of $u_\ep$ given in (\ref{estim_u_ep}), provides
$$\eta_\ep^2R_c\|D w_\ep\|_{L^2(\Omega_\ep)^{3\times 3}}\leq 2N^2\eta_\ep^{5\over 2}C + \eta_\ep^{5\over 2}c_2\|g'\|_{L^2(\omega)^2}.$$
This implies (\ref{estim_w_ep}), and by means of the dilatation, we get (\ref{estim_tilde_w_ep}).\\

Finally, in the case $\eta_\ep\gg \ep$, estimates (\ref{estim_u_ep+})-(\ref{estim_w_ep+}) in $\Omega_\ep^+$ are obtained similarly as above by using the following Poincar\'e's inequality in $\Omega_\ep^+$,
\begin{equation}\label{Poincare+}
\|v\|_{L^2(\Omega_\ep^+)^3}\leq C\ep\|Dv\|_{L^2(\Omega_\ep^+)^{3\times 3}},\quad \forall v\in H^1_0(\Omega_\ep^+)^3.
\end{equation}
This estimate is obtained by using the fact that in the case $\eta_\ep\gg \ep$, in $\Omega_\ep^+$ we can find the boundary with homogeneous boundary condition at distance $\ep$ integrating along the horizontal variable $x'$.

Thus,  taking $u_\ep^+$ as test function in the first equation of (\ref{form_var_1}) and using (\ref{Poincare+}), we get
$$
\begin{array}{l}
\displaystyle
\|Du_\ep^+\|^2_{L^2(\Omega_\ep^+)^{3\times 3}}\\
\noame
\displaystyle \leq 2N^2 \|w_\ep^+\|_{L^2(\Omega_\ep^+)^3}\|D u_\ep^+\|_{L^2(\Omega_\ep^+)^{3\times 3}}+\ep\eta_\ep^{1\over 2}c_2\|f'\|_{L^2(\omega)^2}\|D u_\ep^+\|_{L^2(\Omega_\ep^+)^{3\times 3}}\,,
\end{array}
$$
and then
\begin{equation}\label{lem_estim_6}
\ep^{-1}\eta_\ep^{-{1\over 2}}\|D u_\ep^+\|_{L^2(\Omega_\ep^+)^{3\times 3}}\leq \ep^{-1}\eta_\ep^{-{1\over 2}}2N^2\|w_\ep^+\|_{L^2(\Omega_\ep^+)^3}+c_2\|f'\|_{L^2(\omega)^2}\,.
\end{equation}
Next, we obtain 
\begin{equation*}
\begin{array}{l}
\displaystyle
\eta_\ep^2R_c\|D w_\ep^+\|_{L^2(\Omega_\ep^+)^{3\times 3}}^2+4N^2\|w_\ep^+\|_{L^2(\Omega_\ep^+)^3}^2\\
\noame
\displaystyle\quad\quad 
\leq 2N^2\|w_\ep^+\|_{L^2(\Omega_\ep^+)^3}\|D u_\ep^+\|_{L^2(\Omega_\ep^+)^{3\times 3}}+	\ep\eta_\ep^{1\over 2}\|g'\|_{L^2(\omega)^2}\|w_\ep^+\|_{L^2(\Omega_\ep^+)^3}\,,
\end{array}
\end{equation*}
and then
\begin{equation*}
\ep^{-1}\eta_\ep^{-{1\over 2}}2N^2\|w_\ep^+\|_{L^2(\Omega_\ep^+)^3}\leq \ep^{-1}\eta_\ep^{-{1\over 2}}N^2\|D u_\ep^+\|_{L^2(\Omega_\ep^+)^{3\times 3}}+{1\over 2}\|g'\|_{L^2(\omega)^2}\,.
\end{equation*}
From the above estimates, we get the second estimate in (\ref{estim_u_ep+}) and by (\ref{Poincare+}), the first one. By means of the dilatation we deduce (\ref{estim_tilde_u_ep+}).

Finally, by applying (\ref{Poincare+}), we have
\begin{equation*}
\begin{array}{l}
\displaystyle
\eta_\ep^2R_c\|D w_\ep^+\|_{L^2(\Omega_\ep^+)^{3\times 3}}^2+4N^2\|w_\ep^+\|_{L^2(\Omega_\ep^+)^3}^2\\
\noame
\displaystyle\quad\quad 
\leq 2N^2\|u_\ep^+\|_{L^2(\Omega_\ep^+)^3}\|D w_\ep^+\|_{L^2(\Omega_\ep^+)^{3\times 3}}+\ep^2\eta_\ep^{1\over 2}\|g'\|_{L^2(\omega)^2}\|Dw_\ep^+\|_{L^2(\Omega_\ep^+)^{3\times3}}\,,
\end{array}
\end{equation*}
which, by using the estimate of $u_\ep^+$ given in (\ref{estim_u_ep+}), provides the second estimate in (\ref{estim_w_ep+}), and then the first one. Moreover, by means of the dilatation we deduce (\ref{estim_tilde_w_ep+}) which ends the proof.
\hfill$\square$\end{proof}

\subsection{The extension of $(\tilde u_\ep,\tilde w_\ep,\tilde p_\ep)$ to the whole domain $\Omega$}
The sequence of solutions $(\tilde u_\ep,\tilde w_\ep,\tilde p_\ep)\in H^1_0(\widetilde \Omega_\ep)^3\times H^1_0(\widetilde \Omega_\ep)^3\times L^2_0(\widetilde \Omega_\ep)$ is not defined in a fixed domain independent of $\ep$ but rather in a varying set $\widetilde\Omega_\ep$. In order to pass to the limit if $\ep$ tends to zero, convergences in fixed Sobolev spaces (defined in $\Omega$) are used, which requires first that $(\tilde u_\ep,\tilde w_\ep,\tilde p_\ep)$ be extended to the whole domain $\Omega$.

Therefore, we extend $\tilde{u}_{\varepsilon}$ and $\tilde{w}_{\varepsilon}$  by zero in $\Omega\setminus \widetilde{\Omega}_{\varepsilon}$ (this is compatible with the homogeneous boundary condition on $\partial \widetilde{\Omega}_{\varepsilon}$), and denote the extensions by the same symbol. Obviously, estimates (\ref{estim_u_ep})-(\ref{estim_tilde_w_ep}) remain valid and the extension $\tilde u_\ep$ is divergence free too.

Extending the pressure is a much more difficult task.  A continuation of the pressure for a flow in a porous media was introduced in \cite{Tartar}. This construction applies to periodic holes in a domain $\Omega_\varepsilon$ when each hole is strictly contained into the periodic cell. In this context, we can not use directly this result because the ``holes" are along the boundary $\Sigma_\varepsilon$ of $\Omega_\varepsilon$ and moreover, the scale of the vertical direction is smaller than the scales of the horizontal directions. This fact will induce several limitations in the results obtained by using the method, especially in view of the convergence for the pressure. In this sense, for the case of Newtonian fluids,  an operator $R^\varepsilon$  generalizing the results of \cite{Tartar}  to this context (extending the pressure from $\Omega_\ep$ to $Q_\ep$) was introduced in \cite{Bayada_Chambat,Mikelic2}, and later extended to the case of non-Newtonian (power law) fluids \cite{Anguiano_SG} by defining an extension operator $R^\ep_p$, for every flow index $p>1$. 

Then, in order to extend the pressure to the whole domain $\Omega$, the mapping $R^\ep$  (defined in  \cite[Lemma 4.6]{Anguiano_SG} as $R_2^\ep$) allows us to extend the pressure $p_\ep$ from $\Omega_\ep$ to $Q_\ep$ by introducing $F_\ep$ in $H^{-1}(Q_\ep)^3$ as follows (brackets are for duality products between $H^{-1}$ and $H^1_0$)
\begin{equation}\label{F}\langle F_\varepsilon, \varphi\rangle_{Q_\varepsilon}=\langle \nabla p_\varepsilon, R^\varepsilon (\varphi)\rangle_{\Omega_\varepsilon}\quad \hbox{for any }\varphi\in H^1_0(Q_\varepsilon)^3\,.
\end{equation}
We calcule the right hand side of (\ref{F}) by using the first equation of (\ref{form_var_1}) and we have
\begin{equation}\label{equality_duality}
\begin{array}{rl}
\displaystyle
\left\langle F_{\varepsilon},\varphi\right\rangle_{Q_\varepsilon}=&\displaystyle-\int_{\Omega_\varepsilon}D u_\ep : D R^{\varepsilon}(\varphi)\,dx + 2N^2 \int_{\Omega_\ep}{\rm rot}\,w_\ep\cdot R^\ep(\varphi)\,dx
\\
\noame&\displaystyle
+\int_{\Omega_\varepsilon} f'(x')\cdot R^{\varepsilon}(\varphi)'\,dx\,.
\end{array}\end{equation}
Using Lemma \ref{lemma_estimates} for fixed $\varepsilon$, we see that it is a bounded functional on $H^1_0(Q_\varepsilon)$ (see the proof of Lemma \ref{lemma_est_P} below) and in fact $F_\varepsilon\in H^{-1}(Q_\varepsilon)^3$. Moreover, ${\rm div} \varphi=0$ implies $\left\langle F_{\varepsilon},\varphi\right\rangle_{Q_\varepsilon}=0$, 
and the DeRham theorem gives the existence of $P_\varepsilon$ in $L^{2}_0(Q_\varepsilon)$ with $F_\varepsilon=\nabla P_\varepsilon$.

Defining the rescaled extended pressure $\tilde P_\ep\in L^2_0(\Omega)$ by 
$$\tilde P_\ep(x',y_3)=P_\ep(x',\eta_\ep y_3),\quad\hbox{a.e. }(x',y_3)\in\Omega,$$
we get for any $\tilde \varphi\in H^1_0(\Omega)^3$ where $\tilde\varphi(x',y_3)=\varphi(x',\eta_\ep y_3)$ that 
$$\begin{array}{rl}\displaystyle\langle \nabla_{\eta_\varepsilon}\tilde P_\varepsilon, \tilde\varphi\rangle_{\Omega}&\displaystyle
=-\int_{\Omega}\tilde P_\varepsilon\,{\rm div}_{\eta_\varepsilon}\,\tilde\varphi\,dx'dy_3\\
\noame&\displaystyle
=-\eta_\varepsilon^{-1}\int_{Q_\varepsilon}P_\varepsilon\,{\rm div}\,\varphi\,dx=\eta_\varepsilon^{-1}\langle \nabla P_\varepsilon, \varphi\rangle_{Q_\varepsilon}\,.
\end{array}$$
Then, using the identification (\ref{equality_duality}) of $F_\varepsilon$, we get
$$\begin{array}{rl}\displaystyle\langle \nabla_{\eta_\varepsilon}\tilde P_\varepsilon, \tilde\varphi\rangle_{\Omega}
= &\displaystyle\eta_\varepsilon^{-1}\Big(-\int_{\Omega_\varepsilon}D u_\ep : D R^{\varepsilon}(\varphi)\,dx
+2N^2\int_{\Omega_\ep}{\rm rot}\,w_\ep\cdot R^\ep(\varphi)\,dx
\\
\noame&
\displaystyle
+\int_{\Omega_\varepsilon} f'(x')\cdot R^{\varepsilon}(\varphi)'\,dx\Big)\,,
\end{array}$$
and applying the change of variables (\ref{dilatacion}), we obtain 
\begin{equation}\label{extension_1}
\begin{array}{rl}\displaystyle\langle \nabla_{\eta_\varepsilon}\tilde P_\varepsilon, \tilde\varphi\rangle_{\Omega}
=&\displaystyle- \int_{\widetilde \Omega_\varepsilon}D_{\eta_\ep} \tilde u_\ep : D_{\eta_\ep} \tilde R^{\varepsilon}(\tilde \varphi)\,dx'dy_3
\\
\noame &\displaystyle+2N^2\int_{\widetilde \Omega_\ep}{\rm rot}_{\eta_\ep}\tilde w_\ep\cdot \tilde R^\ep(\tilde \varphi)\,dx'dy_3
\\
\noame&\displaystyle+\int_{\widetilde \Omega_\varepsilon} f(x')\cdot \tilde R^{\varepsilon}(\tilde \varphi)'\,dx'dy_3\,,
\end{array}
\end{equation}
where $\tilde R^\ep(\tilde\varphi)=R^\ep(\varphi)$ for any $ \varphi\in H^1_0(Q_\varepsilon)^3$ where $\tilde\varphi(x',y_3)=\varphi(x',\eta_\ep y_3)$.\\

Now, we estimate the right-hand side of (\ref{extension_1}) to obtain the {\it a priori} estimate of the pressure $\tilde P_\ep$.

\begin{lemma}\label{lemma_est_P}
There exists a constant $C>0$ independent of $\ep$, such that the extension $\tilde P_\ep\in L^2_0(\Omega)$ of the pressure $\tilde p_\ep$ satisfies
\begin{equation}\label{esti_P}
\|\tilde P_\ep\|_{L^2(\Omega)}\leq C.
\end{equation}
\end{lemma}
\begin{proof} From the proof of Lemma 4.7-(i) in \cite{Anguiano_SG}, we have that $\tilde R^\ep(\tilde\varphi)$ satisfies the following estimates
\begin{equation}\label{ext_1}
\begin{array}{l}
\displaystyle\|\tilde R^\ep(\tilde\varphi)\|_{L^2(\widetilde\Omega_\ep)^3}\leq C\left(\|\tilde\varphi\|_{L^2(\Omega)^3}+
\ep\|D_{x'}\tilde\varphi\|_{L^2(\Omega)^{3\times 2}} + \|\partial_{y_3}\tilde\varphi\|_{L^2(\Omega)^3}\right)\,,\\
\noame
\displaystyle
\|D_{x'}\tilde R^\ep(\tilde\varphi)\|_{L^2(\widetilde\Omega_\ep)^{3\times 2}}\leq \! C\!\left({1\over \ep}\|\tilde\varphi\|_{L^2(\Omega)^3}+
\|D_{x'}\tilde\varphi\|_{L^2(\Omega)^{3\times 2}} + {1\over \ep}\|\partial_{y_3}\tilde\varphi\|_{L^2(\Omega)^3}\right)\!,\\
\noame
\displaystyle\|\partial_{y_3}\tilde R^\ep(\tilde\varphi)\|_{L^2(\widetilde\Omega_\ep)^3}\leq C\left(\|\tilde\varphi\|_{L^2(\Omega)^3}+
\ep\|D_{x'}\tilde\varphi\|_{L^2(\Omega)^{3\times 2}} + \|\partial_{y_3}\tilde\varphi\|_{L^2(\Omega)^3}\right)\,.
\end{array}
\end{equation}
Thus, in the cases $\eta_\ep\approx \ep$ or $\eta_\ep\ll\ep$, we have 
\begin{equation}\label{ext_2}
\|\tilde R^\ep(\tilde \varphi)\|_{L^2(\widetilde\Omega_\ep)^3}\leq C\|\tilde \varphi\|_{H^1_0(\Omega)^3},\ \|D_{\eta_\ep}\tilde R^\ep(\tilde \varphi)\|_{L^2(\widetilde\Omega_\ep)^{3\times 3}}\leq {C\over \eta_\ep}\|\tilde \varphi\|_{H^1_0(\Omega)^3}\,,
\end{equation}
and in the case $\eta_\ep\gg\ep$, we have
\begin{equation}\label{ext_3}
\|\tilde R^\ep(\tilde \varphi)\|_{L^2(\widetilde\Omega_\ep)^3}\leq C\|\tilde \varphi\|_{H^1_0(\Omega)^3},\  \|D_{\eta_\ep}\tilde R^\ep(\tilde \varphi)\|_{L^2(\widetilde\Omega_\ep)^{3\times 3}}\leq {C\over \ep}\|\tilde \varphi\|_{H^1_0(\Omega)^3}\,.
\end{equation}
In the cases $\eta_\ep\approx \ep$ or $\eta_\ep\ll\ep$, by using Cauchy-Schwarz's inequality, estimates for $D_{\eta_\ep}\tilde u_\ep$ in (\ref{estim_tilde_u_ep}), for $D_{\eta_\ep}w_\ep$ in (\ref{estim_tilde_w_ep}), $f'\in L^2(\omega)^2$, estimate (\ref{Gaffney}) in $\widetilde\Omega_\ep$, and (\ref{ext_2}), we  obtain
\begin{equation}\label{estim_nabla_p_crit_1}
\begin{array}{rl}
\displaystyle
\left|\int_{\widetilde\Omega_\ep}D_{\eta_\ep}\tilde u_\ep:D_{\eta_\ep}\tilde R^\ep(\tilde\varphi)\,dx'dy_3\right|& \leq \displaystyle C\eta_\ep\|D_{\eta_\ep}\tilde R^\ep(\tilde\varphi)\|_{L^2(\widetilde \Omega_\ep)^{3\times 3}}\\
\noame&\displaystyle\qquad\leq C\|\tilde\varphi\|_{H^1_0(\Omega)^3},\\
\noame
\displaystyle \left|\int_{\widetilde\Omega_\ep}{\rm rot}_{\eta_\ep}w_\ep\cdot \tilde R^\ep(\tilde\varphi)\,dx'dy_3\right| \displaystyle& \leq 
\|D_{\eta_\ep}\tilde w_\ep\|_{L^2(\widetilde\Omega_\ep)^{3\times 3}}\|\tilde R^\ep(\tilde\varphi)\|_{L^2(\widetilde\Omega_\ep)^3}\\
\noame
\displaystyle\qquad  &\leq C\|\tilde R^\ep(\tilde\varphi)\|_{L^2(\widetilde\Omega_\ep)^3}\leq C\|\tilde\varphi\|_{H^1_0(\Omega)^3},\\
\noame\displaystyle
 \left|\int_{\widetilde\Omega_\ep}f'\cdot \tilde R^\ep(\tilde\varphi)\,dx'dy_3\right|&\leq \displaystyle C\|\tilde R^\ep(\tilde\varphi)\|_{L^2(\widetilde\Omega_\ep)^3}\leq C\|\tilde\varphi\|_{H^1_0(\Omega)^3}\,,
\end{array}
\end{equation}
which together with (\ref{extension_1}) gives $\|\nabla_{\eta_\ep}\tilde P_\ep\|_{H^{-1}(\Omega)^3}\leq C$. By using the Ne${\breve{\rm c}}$as inequality there exists a representative $\tilde P_\ep\in L^2_0(\Omega)$ such that
\begin{equation}\label{Necas_inq}
\|\tilde P_\ep\|_{L^2(\Omega)}\leq C\|\nabla\tilde P_\ep\|_{H^{-1}(\Omega)^3}\leq C\|\nabla_{\eta_\ep}\tilde P_\ep\|_{H^{-1}(\Omega)^3},
\end{equation}
which implies (\ref{esti_P}).\\

In the case $\eta_\ep\gg\ep$, due to the highly oscillating boundary, we proceed as the previous cases by considering $\tilde \varphi\in H^1_0(\Omega^+)^3$, estimates (\ref{Gaffney}), (\ref{estim_tilde_u_ep+}) and (\ref{estim_tilde_w_ep+}) in $\widetilde \Omega^+_\ep$ and (\ref{ext_3}), which gives 
$$
\begin{array}{l}
\displaystyle
\left|\int_{\widetilde\Omega_\ep}D_{\eta_\ep}\tilde u^+_\ep:D_{\eta_\ep}\tilde R^\ep(\tilde\varphi)\,dx'dy_3\right|\leq C\ep\|D_{\eta_\ep}\tilde R^\ep(\tilde\varphi)\|_{L^2(\widetilde \Omega_\ep^+)^{3\times 3}}\leq C\|\tilde\varphi\|_{H^1_0(\Omega^+)^3},\\
\noame
\displaystyle \left|\int_{\widetilde\Omega_\ep}{\rm rot}_{\eta_\ep}w_\ep^+\cdot \tilde R^\ep(\tilde\varphi)\,dx'dy_3\right|\leq 
\|D_{\eta_\ep}\tilde w_\ep^+\|_{L^2(\widetilde\Omega_\ep^+)^{3\times 3}}\|\tilde R^\ep(\tilde\varphi)\|_{L^2(\widetilde\Omega_\ep^+)^3}
\\
\noame\displaystyle\qquad \leq C\ep^2\eta_\ep^{-2}\|\tilde R^\ep(\tilde\varphi)\|_{L^2(\widetilde\Omega_\ep^+)^3}\leq C\|\tilde\varphi\|_{H^1_0(\Omega^+)^3},\\
\noame\displaystyle
 \left|\int_{\widetilde\Omega_\ep}f'\cdot \tilde R^\ep(\tilde\varphi)\,dx'dy_3\right|\leq C\|\tilde R^\ep(\tilde\varphi)\|_{L^2(\widetilde\Omega_\ep^+)^3}\leq C\|\tilde\varphi\|_{H^1_0(\Omega^+)^3}\,,
\end{array}
$$
and we deduce $$\|\nabla_{\eta_\ep}\tilde P_\ep\|_{H^{-1}(\Omega^+)^3}\leq C.$$
Finally, reproducing previous computations by considering $\tilde \varphi\in H_0^1(\Omega^-)^3$,  ta\-king into account that $\tilde R^\ep(\tilde\varphi)=\tilde\varphi$ in $\Omega^-$ and estimates (\ref{estim_tilde_u_ep}) and (\ref{estim_tilde_w_ep}) in $\Omega^-$, we deduce that $\|\nabla_{\eta_\ep}\tilde P_\ep\|_{H^{-1}(\Omega^-)^3}\leq C$, which together with the previous estimate, implies $\|\nabla_{\eta_\ep}\tilde P_\ep\|_{H^{-1}(\Omega)^3}\leq C$, and  (\ref{esti_P}) follows from the Ne${\breve{\rm c}}$as inequality (\ref{Necas_inq}).
\hfill$\square$\end{proof}

\subsection{Adaptation of the unfolding method}\label{sec:unfolding}
The change of variables (\ref{dilatacion}) does not provide the information we need about the behavior of $\tilde u_\ep$ and $\tilde w_\ep$ in the microstructure associated to $\widetilde\Omega_\ep$. To solve this difficulty, we use an adaptation of the unfolding method (see \cite{arbogast}, \cite{Ciora}, \cite{Ciora2} for more details) introduced to this context in \cite{Anguiano_SG}.
\\

Let us recall that this adaptation of the unfolding method divides the domain $\widetilde\Omega_\ep$ in cubes of lateral length $\ep$ and vertical length $h(y')$, and the domain $\Omega$ in cubes of lateral length $\ep$ and vertical length $h_{\rm max}$. Thus, given $\tilde{u}_{\varepsilon}, \tilde w_\ep \in H^1_0(\widetilde \Omega_\varepsilon)^3$ the solution of the rescaled system (\ref{system_2})-(\ref{bc_system_2}), we define $\hat{u}_{\varepsilon}$, $\hat{w}_{\varepsilon}$ by
\begin{eqnarray}\label{uhat}
\hat{u}_{\varepsilon}(x^{\prime},y)=\tilde{u}_{\varepsilon}\left( {\varepsilon}\kappa\left(\frac{x^{\prime}}{{\varepsilon}} \right)+{\varepsilon}y^{\prime},y_3 \right)\text{\ \ a.e. \ }(x^{\prime},y)\in \omega\times Y,
\\
\noame
\hat{w}_{\varepsilon}(x^{\prime},y)=\tilde{w}_{\varepsilon}\left( {\varepsilon}\kappa\left(\frac{x^{\prime}}{{\varepsilon}} \right)+{\varepsilon}y^{\prime},y_3 \right)\text{\ \ a.e. \ }(x^{\prime},y)\in \omega\times Y,\label{what}
\end{eqnarray}
and considering the extended pressure $\tilde {P}_{\varepsilon}\in L^2_0(\Omega)$, we define $\hat P_\ep$ by
\begin{eqnarray}
\hat{P}_{\varepsilon}(x^{\prime},y)=\tilde{P}_{\varepsilon}\left( {\varepsilon}\kappa\left(\frac{x^{\prime}}{{\varepsilon}} \right)+{\varepsilon}y^{\prime},y_3 \right)\text{\ \ a.e. \ }(x^{\prime},y)\in \omega\times \Pi,\label{Phat}
\end{eqnarray}
where the functions $\tilde u_\varepsilon$, $\tilde w_\varepsilon$ and $\tilde P_\varepsilon$ are assumed to be  extended by zero outside
$\omega$ and  the function $\kappa$ is defined by (\ref{kappa_fun}).

\begin{remark}\label{remarkCV}
For $k^{\prime}\in T_{\varepsilon}$, the restrictions of $\hat{u}_{\varepsilon}$ and $\hat{w}_{\varepsilon}$  to $Y^{\prime}_{k^{\prime},{\varepsilon}}\times Y$ and $\hat{P}_{\varepsilon}$ to  $Y^{\prime}_{k^{\prime},{\varepsilon}}\times \Pi$ do not depend on $x^{\prime}$, whereas as a function of $y$ it is obtained from $(\tilde{u}_{\varepsilon}, \tilde{P}_{\varepsilon})$ by using the change of variables 
\begin{equation}\label{CV}
y^{\prime}=\frac{x^{\prime}-{\varepsilon}k^{\prime}}{{\varepsilon}}, 
\end{equation}
which transforms $Y_{k^{\prime},{\varepsilon}}$ into $Y$ and  $\widetilde Q_{k',\varepsilon}$ into $\Pi$, respectively.
\end{remark}
We are now in position to obtain estimates for the sequences $(\hat{u}_{\varepsilon}, \hat w_\ep, \hat{P}_{\varepsilon})$. 
 \begin{lemma}\label{estimates_hat}
 There exists a constant $C>0$ independent of $\ep$, such that $\hat u_\ep$, $\hat w_\ep$ and $\hat P_\ep$ defined by (\ref{uhat}), (\ref{what}) and (\ref{Phat}) respectively, satisfy

 \begin{equation}\label{estim_u_hat}
 \begin{array}{c}
 \|\hat u_\ep\|_{L^2(\omega\times Y)^3}\leq C\eta_\ep^2, \quad 
 \|D_{y'}\hat u_\ep\|_{L^2(\omega\times Y)^{3\times 2}}\leq C\ep\eta_\ep,\\
 \noame
 \|\partial_{y_3}\hat u_\ep\|_{L^2(\omega\times Y)^{3}}\leq C\eta_\ep^2,
 \end{array}
 \end{equation}
\begin{equation} \label{estim_w_hat}
 \begin{array}{c}
  \|\hat w_\ep\|_{L^2(\omega\times Y)^3}\leq C\eta_\ep,\quad 
 \|D_{y'}\hat w_\ep\|_{L^2(\omega\times Y)^{3\times 2}}\leq C\ep,\\
 \noame\displaystyle 
 \|\partial_{y_3}\hat w_\ep\|_{L^2(\omega\times Y)^{3}}\leq C\eta_\ep,\\
\end{array}
\end{equation}
\begin{equation}\label{estim_P_hat}
 \|\hat P_\ep\|_{L^2(\omega\times \Pi)^3}\leq C. 
  \end{equation}
 \end{lemma}
\begin{proof}From the proof of Lemma 4.9 in \cite{Anguiano_SG} in the case $p=2$, we have the following properties concerning the estimates of a function $\tilde \varphi_\ep \in H^1_0(\widetilde\Omega_\ep)^3$ and an extended function $\tilde\psi_\ep \in L^2(\Omega)$ and their respective unfolding functions $\hat \varphi_\ep$ and $\hat \psi_\ep$
\begin{eqnarray}
&\|\hat \varphi_\ep\|_{L^2(\omega\times Y)^3}=\|\tilde \varphi_\ep\|_{L^2(\widetilde\Omega_\ep)^3},\quad 
\|D_{y'}\hat \varphi_\ep\|_{L^2(\omega\times Y)^{3\times 2}}=\ep\|D_{x'}\tilde \varphi_\ep\|_{L^2(\widetilde\Omega_\ep)^{3\times 2}}\,,&\nonumber\\
\noame
&\|\hat \psi_\ep\|_{L^2(\omega\times \Pi)}=\|\tilde \psi_\ep\|_{L^2(\Omega)},\quad 
\|\partial_{y_3}\hat \varphi_\ep\|_{L^2(\omega\times Y)^{3}}=\|\partial_{y_3}\tilde \varphi_\ep\|_{L^2(\widetilde\Omega_\ep)^{3}}\,.&\label{relationpressure}
\end{eqnarray}
Thus, combining previous estimates of $\hat \varphi_\ep$ with estimates for $\tilde u_\ep$ and $\tilde w_\ep$ given in (\ref{estim_tilde_u_ep}) and (\ref{estim_tilde_w_ep}), we respectively get (\ref{estim_u_hat}) and (\ref{estim_w_hat}). For the pressure, combining the previous estimate of $\hat \psi_\ep$ with (\ref{esti_P}) we get (\ref{estim_P_hat}).
\hfill$\square$\end{proof}

\paragraph{Weak variational formulation.} To finish this section, we will give the variational formulation satisfied by the functions $(\hat u_\ep,\hat w_\ep,\hat P_\ep)$, which will be useful in the following sections.

We consider $\varphi_\ep(x',y_3)=\varphi(x',x'/\ep,y_3)$ and $\psi_\ep(x',y_3)=\psi(x',x'/\ep,y_3)$ as test function in (\ref{form_var_2}) where $\varphi(x',y)$, $\psi(x',y)\in \mathcal{D}(\omega;C_{\#}^\infty(Y)^3)$, and taking into account the extension of the pressure, we have 
$$\int_{\widetilde\Omega_\ep}\nabla_{\eta_\ep}\tilde p_\ep\cdot \varphi_\ep\,dx'dy_3=\int_{\Omega}\nabla_{\eta_\ep}\tilde P_\ep\cdot \varphi_\ep\,dx'dy_3\,,$$
and so
\begin{equation}\label{form_var_general_1}
 \begin{array}{l}\displaystyle
\int_{\widetilde\Omega_\ep}D_{\eta_\ep}\tilde u_\ep:D_{\eta_\ep}\varphi_\ep\,dx'dy_3-\int_{\Omega}
\tilde P_\ep\,{\rm div}_{\eta_\ep}\varphi_\ep\,dx'dy_3\\
\noame\qquad\displaystyle
=2N^2\int_{\widetilde\Omega_\ep}{\rm rot}_{\eta_\ep}\tilde w_\ep\cdot \varphi_\ep\,dx'dy_3+\int_{\widetilde\Omega_\ep}f'\cdot \varphi_\ep'\,dx'dy_3\,,
\\
\\
\displaystyle
\eta_\ep^2R_c\int_{\widetilde\Omega_\ep}D_{\eta_\ep}\tilde w_\ep:D_{\eta_\ep}\psi_\ep\,dx'dy_3
+4N^2\int_{\widetilde\Omega_\ep}\tilde w_\ep\cdot \psi_\ep\,dx'dy_3\\
\noame
\displaystyle\qquad=2N^2\int_{\widetilde\Omega_\ep}{\rm rot}_{\eta_\ep}\tilde u_\ep\cdot \psi_\ep\,dx'dy_3+\int_{\widetilde\Omega_\ep}g_\ep'\cdot \psi_\ep'\,dx'dy_3\,,
\end{array} 
\end{equation}
where $g_\ep'$ is given by (\ref{estim_f_g_cases_crit_sub}) or (\ref{estim_f_g_cases_sup}) depending on the case.\\

Now, by the change of variables given in Remark \ref{remarkCV} (see \cite{Anguiano_SG} for more details), we obtain
\begin{equation}\label{form_var_hat_u}
 \begin{array}{l}
\displaystyle{1\over \ep^2}\int_{\omega\times Y}D_{y'}\hat u_\ep':D_{y'}\varphi'\,dx'dy+{1\over \eta_\ep^2}\int_{\omega\times Y}\partial_{y_3}\hat u'_\ep :\partial_{y_3}\varphi'\,dx'dy\\
\noame
\displaystyle
-\int_{\omega\times \Pi}\hat P_\ep{\rm div}_{x'}\varphi'\,dx'dy-{1\over \ep}\int_{\omega\times \Pi}\hat P_\ep{\rm div}_{y'}\varphi'\,dx'dy\\
\noame
\displaystyle
={2N^2\over \ep}\int_{\omega\times Y}{\rm rot}_{y'}\hat w_{\ep,3}\cdot \varphi'\,dx'dy+
{2N^2\over \eta_\ep}\int_{\omega\times Y}{\rm rot}_{y_3}\hat w_\ep'\cdot \varphi'\,dx'dy\\
\noame\displaystyle+
\int_{\omega\times Y}f'\cdot \varphi'\,dx'dy+O_\ep\,,
\\
\\
\displaystyle
\displaystyle{1\over \ep^2}\int_{\omega\times Y}\nabla_{y'}\hat u_{\ep,3}\cdot\nabla_{y'}\varphi_3\,dx'dy+{1\over \eta_\ep^2}\int_{\omega\times Y}\partial_{y_3}\hat u_{\ep,3} \cdot \partial_{y_3}\varphi_3\,dx'dy\\
\noame
\displaystyle
-{1\over \eta_\ep}\int_{\omega\times \Pi}\hat P_\ep\partial_{y_3}\varphi_3\,dx'dy
={2N^2\over \ep}\int_{\omega\times Y}{\rm Rot}_{y'}\hat w_{\ep}'\, \varphi_3\,dx'dy+O_\ep\,,
\end{array}
\end{equation}
and
\begin{equation}\label{form_var_hat_w}
 \begin{array}{l}
\displaystyle{\eta_\ep^2\over \ep^2}R_c\int_{\omega\times Y}D_{y'}\hat w_\ep':D_{y'}\psi'\,dx'dy+R_c\int_{\omega\times Y}\partial_{y_3}\hat w'_\ep :\partial_{y_3}\psi'\,dx'dy
\\
\noame\displaystyle+ 4N^2\int_{\omega\times Y}\hat w_\ep'\cdot \psi'\,dx'dy\\
\noame
\displaystyle
={2N^2\over \ep}\int_{\omega\times Y}{\rm rot}_{y'}\hat u_{\ep,3}\cdot \psi'\,dx'dy+
{2N^2\over \eta_\ep}\int_{\omega\times Y}{\rm rot}_{y_3}\hat u_\ep'\cdot \psi'\,dx'dy
\\
\noame\displaystyle
+
\int_{\omega\times Y}g_\ep'\cdot \psi'\,dx'dy+O_\ep\,,
\\
\\
\displaystyle
\displaystyle{\eta_\ep^2\over \ep^2}R_c\int_{\omega\times Y}\nabla_{y'}\hat w_{\ep,3}\cdot \nabla_{y'}\psi_3\,dx'dy+R_c\int_{\omega\times Y}\partial_{y_3}\hat w_{\ep,3} :\partial_{y_3}\psi_3\,dx'dy\\
\noame
\displaystyle
+4N^2\int_{\omega\times Y}\hat w_{\ep,3}\cdot \psi_3\,dx'dy={2N^2\over \ep}\int_{\omega\times Y}{\rm Rot}_{y'}\hat u_{\ep}'\, \psi_3\,dx'dy+O_\ep\,.
\end{array}
\end{equation}

When $\ep$ tends to zero, we obtain for $(\hat u_\ep, \hat w_\ep, \hat P_\ep)$ different asymptotic behaviors, depending on the magnitude of $\eta_\ep$ with respect to $\ep$. We will analyze them in the next sections.

\section{Stokes roughness regime ($0<\lambda<+\infty$)}\label{SRR}
It corresponds to the critical case when the thickness of the domain is proportional to the wavelength of the roughness, with $\lambda$ the proportionality constant, that is $\eta_\ep\approx\ep$, with $\eta_\ep/\ep\to \lambda$, $0<\lambda<+\infty$.

Let us introduce some notation which will be useful along this section. For a vectorial function $v=(v',v_3)$ and a scalar function $w$, we introduce the operators $D_\lambda$, $\nabla_\lambda$, ${\rm div}_\lambda$ and ${\rm rot}_\lambda$ by 
\begin{eqnarray}
&(D_{\lambda}v)_{ij}=\lambda \partial_{x_j}v_i\hbox{ for }i=1,2,3,\ j=1,2,\quad (D_{\lambda}v)_{i,3}=\partial_{y_3}v_i\hbox{ for }i=1,2,3,\nonumber&\\
\noame
&\Delta_\lambda v=\lambda^2\Delta_{y'}v+\partial_{y_3}^2 v,\quad \nabla_{\lambda} w=(\lambda\nabla_{y'}w,\partial_{y_3} w)^t,&\nonumber\\
\noame
&{\rm div}_\lambda v=\lambda{\rm div}_{y'}v'+\partial_{y_3} v_3,\quad {\rm rot}_\lambda v=(\lambda{\rm rot}_{y'}v_3+
{\rm rot}_{y_3} v',\lambda {\rm Rot}_{y'} v')\,,&\nonumber
\end{eqnarray}
where ${\rm rot}_{y'}$, ${\rm rot}_{y_3}$ and ${\rm Rot}_{y'}$ are defined in (\ref{def_rot}).
Next, we give some compactness results about the behavior of the extended sequences $(\tilde u_\ep, \tilde w_\ep, \tilde P_\ep)$ and the related unfolding functions $(\hat u_\ep, \hat w_\ep, \hat P_\ep)$ satisfying the {\it a priori} estimates given in Lemmas \ref{lemma_estimates} and \ref{lemma_est_P}, and Lemma \ref{estimates_hat} respectively.
\begin{lemma}\label{lem_asymp_crit}
For a subsequence of $\ep$ still denote by $\ep$, we have that 
\begin{itemize}
\item[(i)] {\it (Velocity)} there exist $\tilde u\in H^1(0,h_{\rm max};L^2(\omega)^3)$, with $\tilde u=0$ on $y_3=\{0,h_{\rm max}\}$ and $\tilde u_3=0$, and  $\hat u\in L^2(\omega;H^1_{\#}(Y))^3$, with $\hat u=0$ on $y_3=\{0,h(y')\}$ such that $\int_{Y}\hat u(x',y)dy=\int_0^{h_{\rm max}}\tilde u(x',y_3)\,dy_3$ with $\int_{Y}\hat u_3\,dy=0$, and moreover 
\begin{equation}\label{conv_u_crit}
\begin{array}{c}
\displaystyle\eta_\ep^{-2}\tilde u_\ep\rightharpoonup (\tilde u',0)\hbox{  in  }H^1(0,h_{\rm max};L^2(\omega)^3),\\
\noame
\eta_\ep^{-2}\hat u_\ep\rightharpoonup \hat u\hbox{  in  }L^2(\omega;H^1(Y)^3),
\end{array}
\end{equation}
\begin{equation}\label{div_x_crit}
\begin{array}{c}
\displaystyle{\rm div}_{x'}\left(\int_0^{h_{\rm max}}\tilde u'(x',y_3)\,dy_3\right)=0\hbox{  in  }\omega,\\
\noame\displaystyle
\left(\int_0^{h_{\rm max}}\tilde u'(x',y_3)\,dy_3\right)\cdot n=0\hbox{  in  }\partial\omega\,,
\end{array}
\end{equation}
\begin{equation}\label{div_y_crit}
\begin{array}{c}
\displaystyle{\rm div}_\lambda\hat u=0\hbox{  in  }\omega\times Y,\\
\noame
\displaystyle
{\rm div}_{x'}\left(\int_{Y}\hat u'(x',y)\,dy\right)=0\hbox{  in  }\omega,\left(\int_{Y}\hat u'(x',y)\,dy\right)\cdot n=0
\hbox{  in  }\partial\omega\,,
\end{array}
\end{equation}
\item[(ii)] {\it (Microrotation)} there exist $\tilde w\in H^1(0,h_{\rm max};L^2(\omega)^3)$, with $\tilde w=0$ on $y_3=\{0,h_{\rm max}\}$ and $\tilde w_3=0$, and  $\hat w\in L^2(\omega;H^1_{\#}(Y))^3$, with $\hat w=0$ on $y_3=\{0,h(y')\}$ such that $\int_{Y}\hat w(x',y)dy=\int_0^{h_{\rm max}}\tilde w(x',y_3)\,dy_3$ with $\int_{Y}\hat w_3\,dy=0$, and moreover 
\begin{equation}\label{conv_w_crit}
\begin{array}{c}
\displaystyle\eta_\ep^{-1}\tilde w_\ep\rightharpoonup (\tilde w',0)\hbox{  in  }H^1(0,h_{\rm max};L^2(\omega)^3),\\
\noame
\displaystyle
\eta_\ep^{-1}\hat w_\ep\rightharpoonup \hat w\hbox{  in  }L^2(\omega;H^1(Y)^3),
\end{array}
\end{equation}
\item[(iii)] {\it (Pressure)} there exists a function $\tilde P\in L^2_0(\Omega)$, independent of $y_3$, such that
\begin{eqnarray}
&\displaystyle\tilde P_\ep\to \tilde P\hbox{  in  }L^2(\Omega),\quad \hat P_\ep\to \tilde P\hbox{  in  }L^2(\omega\times \Pi).& \label{conv_P_crit}
\end{eqnarray}
\end{itemize}
\end{lemma}
\begin{proof} We start proving $(i)$. We will only give some remarks and,  for more details,  we refer the reader to Lemmas 5.2-i) and 5.4-i) in \cite{Anguiano_SG}.

We start with the extension $\tilde u_\ep$. Estimates  (\ref{estim_tilde_u_ep}) imply the existence of $\tilde u\in H^1(0,h_{\rm max};L^2(\omega)^3)$ such that convergence (\ref{conv_u_crit})$_1$ holds, and the continuity of the trace applications from the space of $\tilde u$ such that $\|\tilde u\|_{L^2}$ and $\|\partial_{y_3}\tilde u\|_{L^2}$ are bounded to $L^2(\Sigma)$ and to $L^2(\omega\times \{0\})$ implies $\tilde u=0$ on $\Sigma$ and $\omega\times \{0\}$. Next, from the free divergence condition ${\rm div}_{\eta_\ep}\tilde u_\ep=0$, it can be deduced that $\tilde u_3$ is independent of $y_3$, which together with the boundary conditions satisfied by $\tilde u_3$ on $y_3=\{0,h_{\rm max}\}$ implies that $\tilde u_3=0$. Finally, from the free divergence condition and the convergence (\ref{conv_u_crit})$_1$ of $\tilde u_\ep$, it is straightforward the corresponding free divergence condition in a thin domain given in (\ref{div_x_crit}).\\
 
Concerning $\hat u_\ep$, estimates given in (\ref{estim_u_hat}) imply the existence of a function $\hat u\in  L^2(\omega;H^1(Y)^3)$ such that convergence (\ref{conv_u_crit})$_2$ holds. It can be proved the $Y'$-periodicity of $\hat u$, and applying the change of variables (\ref{CV}) to the free divergence condition ${\rm div}_{\eta_\ep}\tilde u_\ep=0$, passing to the limit and taking into account that $\eta_\ep/\ep\to \lambda$, we get divergence condition ${\rm div}_\lambda \hat u=0$ given in (\ref{div_y_crit}). Finally, 
 it can be proved that $\int_Y\hat u(x',y)\,dy=\int_0^{h_{\rm max}}\tilde u(x',y_3)\,dy_3$ which together with $\tilde u_3=0$ implies $\int_0^{h_{\rm max}}\tilde u_3(x',y_3)\,dy_3=0$, and together with property (\ref{div_x_crit}) implies the divergence condition ${\rm div}_{x'}\int_Y\hat u'(x',y)dy=0$ given in (\ref{div_y_crit}).\\

We continue proving $(ii)$.  From estimates (\ref{estim_tilde_w_ep}), the first convergence of (\ref{conv_w_crit}) and that $\tilde w=0$ on $y_3=\{0,h_{\rm max}\}$ straighfordward. It remains to prove that $\tilde w_3=0$. To do this, we consider as test function $\psi_\ep(x',y_3)=(0,0,\eta_\ep^{-1}\psi_3)$ in the variational formulation (\ref{form_var_general_1}) extended to $\Omega$, and we get
$$
\begin{array}{l}
\displaystyle\eta_\ep R_c\int_{\Omega}\nabla_{x'}\tilde w_{\ep,3}\cdot \nabla_{x'}\psi_3\,dx'dy_3+\eta_\ep^{-1}R_c\int_\Omega
\partial_{y_3}\tilde w_{\ep,3}\,\partial_{y_3}\psi_3\,dx'dy_3\\
\noame
\displaystyle+4N^2\eta_\ep^{-1}\int_\Omega\tilde w_{\ep,3}\psi_3\,dx'dy_3=\eta_\ep^{-1}\int_{\Omega}{\rm Rot}_{x'}\tilde u'_\ep\,\psi_3\,dx'dy_3\,.
\end{array}
$$
Passing to the limit by using concergences of $\tilde u_\ep$ and $\tilde w_\ep$ given in (\ref{conv_u_crit}) and (\ref{conv_w_crit}), we get
$$R_c\int_\Omega \partial_{y_3}\tilde w_3\,\partial_{y_3}\psi_3\,dx'dy_3+4N^2\int_\Omega \tilde w_3\,\psi_3\,dx'dy_3=0\,,$$
and taking into account that $\tilde w_3=0$ on $y_3=\{0,h_{\rm max}\}$,  it is easily deduced that $\tilde w_3=0$ a.e. in $\Omega$.

The proofs of the convergence of $\hat w_\ep$ and identity $\int_Y\hat u\,dy=\int_0^{h_{\rm max}}\tilde w\,dy_3$ are similar to the ones of $\hat u_\ep$ just taking into account estimate (\ref{estim_w_hat}).

We finish the proof with $(iii)$. Estimate (\ref{estim_P_hat}) implies, up to a subsequence, the existence of $\tilde P\in L^2_0(\Omega)$ such that 
\begin{equation}\label{conv_p_weak}
\tilde P_\ep\rightharpoonup \tilde P\quad \hbox{  in  }L^2(\Omega).
\end{equation} Also, from $\|\nabla_{\eta_\ep}\tilde P_\ep\|_{H^{-1}(\Omega)^3}\leq C$, by noting that $\partial_{y_3}\tilde P_\ep/\eta_\ep$ also converges weakly in $H^{-1}(\Omega)$, we obtain $\partial_{y_3}\tilde P=0$ and so $\tilde P$ is independent of $y_3$.\\

 Next, following  \cite{Tartar}, we prove that the convergence of the pressure is in fact strong.  Let $\sigma_\ep\in H^1_0(\Omega)^3$ be such that $
\sigma_\ep\rightharpoonup \sigma\quad\hbox{in }H^1_0(\Omega)^3
$.  Denoting $\tilde\sigma_\ep=(\sigma'_\ep,\ep \sigma_{\ep,3})$ and $\tilde\sigma=(\sigma',0)$, we have
\begin{equation}\label{strong_p_1}
\tilde \sigma_\ep\rightharpoonup \tilde \sigma\quad\hbox{in }H^1_0(\Omega)^3, 
\end{equation}
Then, as $\tilde P$ only depends on $x'$ and denoting $\nabla_{x',y_3}=(\nabla_{x'},\partial_{y_3})^t$, we have
$$\begin{array}{l}
\displaystyle
\left|<\nabla_{x',y_3}\tilde P_\ep,\sigma_\ep>_{\Omega}-<\nabla_{x'}\tilde P,\tilde \sigma>_{\Omega}\right| 
\\
\noame
 \displaystyle \leq 
\left|<\nabla_{x',y_3}\tilde P_\ep-\nabla_{x'}\tilde P,\tilde \sigma>_{\Omega}\right|+\left|<\nabla_{x',y_3}\tilde P_\ep,\sigma_\ep-\tilde \sigma>_{\Omega}\right|.
\end{array}$$
On the one hand, using convergence (\ref{conv_p_weak}), we have 
$$\left|<\nabla_{x',y_3} \tilde P_\ep-\nabla_{x'}\tilde P,\tilde \sigma>_{\Omega}\right|=\left|\int_\Omega\left(\tilde P_\ep-\tilde P\right)\,{\rm div}_{x'}\sigma'\,dx\right|\to 0,\quad \hbox{as }\ep\to 0\,.$$
On the other hand, proceeding as in the proof of Lemma \ref{lemma_est_P}, we have
$$\begin{array}{l}
\left|<\nabla_{x',y_3}\tilde P_\ep,\sigma_\ep-\tilde\sigma>_{\Omega}\right|= \left|<\nabla_{\eta_\ep}\tilde P_\ep,\tilde\sigma_\ep-\tilde\sigma>_{\Omega}\right|\\
\noame\displaystyle
 \leq C\left(\|\tilde \sigma_\ep-\tilde \sigma\|_{L^2(\Omega)^3}+\eta_\ep\|D_{x',y_3} (\sigma_\ep-\tilde \sigma)\|_{L^2(\Omega)^{3\times 3}} \right)\to 0,\quad \hbox{as }\ep\to 0\,,$$
\end{array}$$
 by virtue of $\eta_\varepsilon$ tends to zero,  (\ref{strong_p_1}) and the Rellich theorem. This implies that $\nabla_{x',y_3}\tilde P_\ep\to \nabla_{x'} \tilde P$ strongly in $H^{-1}(\Omega)^3$, which together the classical Ne${\breve{\rm c}}$as inequality implies the strong convergence of the pressure $\tilde P_\ep$ given in (\ref{conv_P_crit}).  Finally, we remark that the strong convergence of sequence $\hat P_\ep$ to $\tilde P$ is a consequence of the strong convergence of $\tilde P_\ep$ to $\tilde P$ (see  \cite[Proposition 2.9]{Ciora2}). 
\hfill$\square$\end{proof}

Using previous convergences, in the following theorem we give the homo\-genized system satisfied by $(\hat u, \hat w, \tilde P)$.
\begin{theorem}
In the case $\eta_\ep\approx \ep$, with $\eta_\ep/\ep\to \lambda$, $0<\lambda<+\infty$, then the sequence $(\eta_\ep^{-2}\hat u_\ep, \eta_\ep^{-1}\hat w_\ep)$ converges weakly to $(\hat u,\hat w)$ in $L^2(\omega;H^1(Y)^3)\times L^2(\omega;$ $H^1(Y)^3)$ and  $\hat P_\ep$ converges strongly to $\tilde P$ in $L^2(\Omega)$, where $(\hat u,\hat w, \tilde P)\in L^2(\omega;$ $H^1_{\#}(Y)^3)\times L^2(\omega;H^1_{\#}(Y)^3)\times (L^2_0(\omega)\cap H^1(\omega))$, with $\int_Y\hat u_3\,dy=\int_Y\hat w_3\,dy=0$, is the unique solution of the following homogenized system
\begin{equation}\label{hom_system_crit}
\left\{\begin{array}{rl}
\displaystyle
-\Delta_\lambda \hat u+\nabla_\lambda \hat q=2N^2 {\rm rot}_\lambda\hat w+f'(x')-\nabla_{x'}\tilde P(x')&\hbox{ in }\omega\times Y,\\
\noame
{\rm div}_\lambda\hat u=0&\hbox{ in }\omega\times Y,\\
\noame
-R_c\Delta_\lambda \hat w+4N^2 \hat w=2N^2 {\rm rot}_\lambda\hat u+g'(x')&\hbox{ in }\omega\times Y,\\
\noame
\hat u=0&\hbox{ on }y_3=\{0,h(y')\},\\
\noame
\displaystyle{\rm div}_{x'}\left(\int_{Y}\hat u'(x',y)\,dy\right)=0&\hbox{ in }\omega,\\
\noame
\displaystyle\left(\int_{Y}\hat u'(x',y)\,dy\right)\cdot n=0&\hbox{ on }\partial\omega,\\
\noame
 \hat q(x',y)\in L^2(\omega;L^2_{0,\#}(Y)).&
\end{array}\right.
\end{equation}
\end{theorem}
\begin{proof} From Lemma \ref{lem_asymp_crit}, conditions (\ref{hom_system_crit})$_{2,4,5,6}$ hold. To prove that $(\hat u, \hat w, \tilde P)$ satisfies the momentum equations given in (\ref{hom_system_crit}), we consider  $\varphi\in \mathcal{D}(\omega;C_{\#}^\infty(Y)^3)$ with ${\rm div}_\lambda\varphi=0$ in $\omega\times Y$ and ${\rm div}_{x'}(\int_Y \varphi'\,dy)=0$ in $\omega$, and we choose $\varphi_\ep=(\lambda(\ep/\eta_\ep)\varphi',\varphi_3)$ in (\ref{form_var_hat_u}). Taking into account that thanks to ${\rm div}_\lambda\varphi=0$ in $\omega\times Y$, we have that 
$${1\over \eta_\ep}\int_{\omega\times \Pi} \hat P_\ep (\lambda{\rm div}_{y'}\varphi'+\partial_{y_3}\varphi_3)\,dx'dy=0\,.$$
Thus, passing to the limit using the convergences (\ref{conv_u_crit}) and (\ref{conv_w_crit}), and taking into account that $\lambda(\ep/\eta_\ep)\to 1$, we obtain
\begin{equation}\label{form_varl_limit_crit_con_p}
\begin{array}{l}
\displaystyle\int_{\omega\times Y}D_\lambda\hat u:D_y \hat \varphi\,dx'dy-\int_{\omega\times \Pi}\tilde P\,{\rm div}_{x'}\varphi'\,dx'dy\\
\noame
\displaystyle
=2N^2 \int_{\omega\times Y}\left(\lambda{\rm rot}_{y'}\hat w_3\cdot \varphi'+{\rm rot}_{y_3}\hat w'\cdot \varphi'+\lambda{\rm rot}_{y'}\hat w'\,\varphi_3\right)dx'dy\\
\noame
\displaystyle+\int_{\omega\times Y}f'\cdot \varphi'\,dx'dy\,.
\end{array}
\end{equation}
Since $\tilde P$ does not depend on $y$ and ${\rm div}_{x'}\int_{Y}\varphi'\,dy=0$ in $\omega$, we have that
$$\int_{\omega\times Y}\tilde P\,{\rm div}_{x'}\varphi'\,dx'dy=\int_\omega \tilde P\,{\rm div}_{x'}\left(\int_{Y}\varphi'\,dy\right)dx'=0,$$
so we get
\begin{equation}\label{form_var_limit_crit}\int_{\omega\times Y}D_\lambda\hat u:D_y\varphi\,dx'dy=2N^2\int_{\omega\times Y}{\rm rot}_{\lambda}\hat w\cdot \varphi\,dx'dy+\int_{\omega\times Y}f'\cdot \varphi'\,dx'dy\,.
\end{equation}
Next, for every $\psi\in \mathcal{D}(\omega;C_\#^\infty(Y)^3)$, we choose $\psi_\ep=\eta_\ep^{-1}\psi$ in (\ref{form_var_hat_w}). Then, passing to the limit using convergences (\ref{conv_u_crit}) and (\ref{conv_w_crit}), we get
\begin{equation}\label{form_var_limit_crit_w}
\begin{array}{l}\displaystyle\int_{\omega\times Y}D_\lambda\hat w:D_y\psi\,dx'dy+4N^2\int_{\omega\times Y}\hat w\cdot \psi\,dx'dy 
\\ \noame
\displaystyle
=2N^2\int_{\omega\times Y}{\rm rot}_{\lambda}\hat u\cdot \psi\,dx'dy+\int_{\omega\times Y}g'\cdot \psi'\,dx'dy\,.
\end{array}
\end{equation}
By density (\ref{form_var_limit_crit}) holds for every function $\varphi$ in the Hilbert space $V$ defined by
$$
V=\left\{\begin{array}{l}
\varphi(x',y)\in L^2(\omega;H^1_{\#}(Y)^3), \hbox{ such that }{\rm div}_\lambda\varphi(x',y)=0\hbox{ in }\omega\times Y,\\
\noame
\displaystyle {\rm div}_{x'}\left(\int_{Y}\varphi'(x',y)\,dy \right)=0\hbox{ in }\omega,\  \left(\int_{Y}\varphi'(x',y)\,dy \right)\cdot n=0\hbox{ on }\partial\omega
\end{array}\right\},
$$
and (\ref{form_var_limit_crit_w}) in $L^2(\omega;H^1_{\#}(Y)^3)$.

From  \cite[Part III, Theorem 2.4.2]{Luka}, the variational formulation (\ref{form_var_limit_crit})-(\ref{form_var_limit_crit_w}) admits a unique solution $(\hat u, \hat w)$ in $V\times L^2(\omega;H^1_{\#}(Y)^3)$. Following  \cite{Allaire0}, the orthogonal of $V$ with respect to the usual scalar product in $L^2(\omega\times Y)$ is made of gradients of the form $\nabla_{x'}q(x')+\nabla_{\lambda}\hat q(x',y)$, with $q(x')\in L^2(\omega)/\mathbb{R}$ and $\hat q(x',y)\in L^2(\omega;L^2_\#(Y)/\mathbb{R})$. Therefore, by integration by parts, the variational formulations (\ref{form_var_limit_crit})-(\ref{form_var_limit_crit_w})  are equivalent to the homogenized system (\ref{hom_system_crit}). It remains to prove that the pressure $q(x')$, arising as a Lagrange multiplier of the incompressibility constraint ${\rm div}_{x'}(\int_Y \hat u'(x',y)dy)=0$, is the same as the limit of the pressure $\hat P_\ep$. This can be easily done by considering in equation (\ref{form_var_hat_u}) a test function only with ${\rm div}_\lambda$ equal to zero,  obtain  the variational formulation (\ref{form_varl_limit_crit_con_p}) and indentifying limits. Since $2N^2{\rm rot}_\lambda\hat w+ f'\in L^2(\omega\times Y)^3$ and $Y$ is smooth enough, we deduce that $\tilde P\in H^1(\omega)$.

Finally, since from  \cite[Part III, Lemma 2.4.1]{Luka}  we have that  (\ref{hom_system_crit}) admits a unique solution, and then the complete sequence $(\eta_\ep^2\hat u_\ep, \eta_\ep^{-1}\hat w_\ep, \hat P_\ep)$ converges to the unique solution $(\hat u(x',y), \hat w(x',y), \tilde P(x'))$. 
\hfill$\square$\end{proof}

Let us define the local problems which are useful to eliminate the variable $y$ of the previous homogenized problem and then obtain a Reynolds equation for the pressure $\tilde P$.

For every $i,k= 1, 2$ and $0<\lambda<+\infty$, we consider the following 3D local micropolar problems  
\begin{equation}\label{local_problems_crit}
\left\{\begin{array}{rl}\displaystyle
-\Delta_\lambda u^{i,k}+\nabla_{\lambda}\pi^{i,k}-2N^2{\rm rot}_\lambda w^{i,k}=e_i\delta_{1k}&\hbox{ in }Y,\\
\noame
{\rm div}_\lambda u^{i,k}=0&\hbox{ in }Y,\\
\noame
-R_c\Delta_\lambda w^{i,k}+4N^2 w^{i,k}-2N^2{\rm rot}_\lambda u^{i,k}=e_i\delta_{2k}&\hbox{ in }Y,\\
\noame
u^{i,k}=w^{i,k}=0&\hbox{ on }y_3=\{0,h(y')\},\\
\noame
u^{i,k}(y), w^{i,k}(y),\pi^{i,k}(y)\quad Y'-\hbox{periodic}.
\end{array}\right.
\end{equation}
It is known (see  \cite[Part III, Lemma 2.5.1]{Luka}) that there exist a unique solution $(u^{i,k},w^{i,k},\pi^{i,k})\in H^1_\#(Y)^3\times H^1_\#(Y)^3\times L^2_0(Y)$ of problem (\ref{local_problems_crit}), and moreover $\pi^{i,k}\in H^1(Y)$.

We give the main result concerning the homogenized flow.  

\begin{theorem}\label{them_main_crit}
Let $(\hat u,\hat w,\tilde P)\in L^2(\omega;H^1_\#(Y)^3)\times L^2(\omega;H^1_\#(Y)^3)\times (L_0^2(\omega)\cap H^1(\omega))$ be the unique weak solution of problem (\ref{hom_system_crit}). Then, the extensions $(\eta_\ep^{-2}\tilde u_\ep,$ $\eta_\ep^{-1}\tilde w_\ep)$ and  $\tilde P_\ep$ of the solution of problem (\ref{system_2})-(\ref{bc_system_2}) converge weakly to $(\tilde u,\tilde w)$   in $H^1(\omega,h_{\rm max};L^2(\omega)^3)\times H^1(\omega,h_{\rm max};L^2(\omega)^3)$ and strongly to $\tilde P$ in $L^2(\Omega)$ res\-pectively, with $\tilde u_3=\tilde w_3=0$. Moreover, defining $\widetilde U(x')=\int_0^{h_{\rm max}}\tilde u(x',y_3)\,dy_3$ and $\widetilde W(x')=\int_0^{h_{\rm max}}\tilde w(x',y_3)\,dy_3$, it holds 
\begin{equation}\label{Darcy_law_u_w_crit}
\begin{array}{ll}
\widetilde U'(x')=K^{(1)}_\lambda\left(f'(x')-\nabla_{x'}\tilde P(x')\right)+K^{(2)}_\lambda g(x'),&\quad \widetilde U_3(x')=0\quad \hbox{ in }\omega,\\
\noame
\widetilde W'(x')=L^{(1)}_\lambda\left(f'(x')-\nabla_{x'}\tilde P(x')\right)+L^{(2)}_\lambda g(x'),&\quad \widetilde W_3(x')=0\quad \hbox{ in }\omega,
\end{array}
\end{equation}
where $K^{(k)}_\lambda$, $L^{(k)}_\lambda\in \mathbb{R}^{2\times 2}$, $k=1,2$, are matrices with coefficients
$$ \left(K^{(k)}_\lambda\right)_{ij}=\int_Y u^{i,k}_j(y)\,dy,\quad  \left(L^{(k)}_\lambda\right)_{ij}=\int_Y w^{i,k}_j(y)\,dy,\quad i,j=1,2,$$
with $u^{i,k}$, $w^{i,k}$, $i,k=1,2$, the solutions of the local micropolar problems defined in (\ref{local_problems_crit}). 

Here, $\tilde P\in H^1(\omega)\cap L^2_0(\omega)$ is the unique solution of the Reynolds problem
\begin{equation}\label{Darcy_law_P_crit}
\left\{\begin{array}{l}
{\rm div}_{x'}\left( -A_\lambda \nabla_{x'}\tilde P(x')+ b_\lambda(x')\right)=0\quad \hbox{ in }\omega,\\
\noame
\left( -A_\lambda \nabla_{x'}\tilde P(x')+ b_\lambda(x')\right)\cdot n=0\quad \hbox{ on }\partial\omega,
\end{array}\right.
\end{equation}
where the flow factors are given by $A_\lambda=K_\lambda^{(1)}$ and $b_\lambda(x')=K_\lambda^{(1)}f'(x')+K_\lambda^{(2)}g'(x')$.
\end{theorem}
\begin{proof}We eliminate the microscopic variable $y$ in the effective problem (\ref{hom_system_crit}). To do that, we consider the following identification
\begin{eqnarray}
&&\hat u'(x',y)=\sum_{i=1}^2\left[\left(f_i(x')-\partial_{x_i}\tilde P(x')\right) (u^{i,1})'(y)+ g_i(x') (u^{i,2})'(y)\right],\nonumber\\
\noame
&&\displaystyle \hat w'(x',y)=\sum_{i=1}^2\left[\left(f_i(x')-\partial_{x_i}\tilde P(x')\right) (w^{i,1})'(y)+ g_i(x') (w^{i,2})'(y)\right],\nonumber\\
\noame
&&\hat q(x',y)=\sum_{i=1}^2\left[\left(f_i(x')-\partial_{x_i}\tilde P(x')\right) \pi^{i,1}(y) + g_i(x') \pi^{i,2}(y)\right].\nonumber
\end{eqnarray}
From to the identities for the  velocity $\int_Y \hat u'(x',y)\,dy=\int_0^{h_{\rm max}}\tilde u'(x',y_3)\,dy_3$ and $\int_Y\hat u_3\,dy=0$, and for the microrotation $\int_Y\hat w'(x',y)\,dy$ $=\int_0^{h_{\rm max}}\tilde w' (x',y_3)\,dy_3$ and  
$\int_Y\hat w_3\,dy=0$ given in Lemma \ref{lem_asymp_crit}, we deduce that $\widetilde U$ and $\widetilde W$ are given by (\ref{Darcy_law_u_w_crit}).

Finally, the divergence condition with respect to the variable $x'$ given in  (\ref{hom_system_crit}) together with the expression of $\widetilde U'(x')$ gives (\ref{Darcy_law_P_crit}).
\hfill$\square$\end{proof}

\section{Reynolds roughness regime ($\lambda=0$)}\label{RRR}
It corresponds to the case when the wavelength of the roughness is much greater than the film thickness, i.e. $\eta_\ep\ll \ep$ which is equivalent to $\lambda=0$.

Next, we give some compactness results about the behavior of the extended sequences $(\tilde u_\ep,\tilde w_\ep,\tilde P_\ep)$ and the unfolding functions $(\hat u_\ep,\hat w_\ep, \hat P_\ep)$ satisfying the {\it a priori} estimates given in Lemmas \ref{lemma_estimates} and \ref{lemma_est_P}, and Lemma \ref{estimates_hat}, respectively.

\begin{lemma}\label{lem_asymp_sub}
For a subsequence of $\ep$ still denoted by $\ep$, we have that
\begin{itemize}
\item[(i)] (Velocity) there exist $\tilde u\in H^1(0,h_{\rm max};L^2(\omega)^3)$, with $\tilde u=0$ on $y_3=\{0,h_{\rm max}\}$ and $\tilde u_3=0$, and  $\hat u\in H^1(0,h(y'); L^2_{\#}(\omega\times Y')^3)$, with  $\hat u=0$ on $y_3=\{0,h(y')\}$ and $\hat u_3$ independent of $y_3$, such that $\int_{Y}\hat u(x',y)dy=$\linebreak$\int_0^{h_{\rm max}}\tilde u(x',y_3)\,dy_3$ with $\int_{Y}\hat u_3\,dy=0$, and moreover
\begin{equation}\label{conv_u_sub}
\begin{array}{c}
\displaystyle\eta_\ep^{-2}\tilde u_\ep\rightharpoonup (\tilde u',0)\hbox{  in  }H^1(0,h_{\rm max};L^2(\omega)^3),\\ 
\noame
\displaystyle
\eta_\ep^{-2}\hat u_\ep\rightharpoonup \hat u\hbox{  in  } H^1(0,h(y'); L^2(\omega\times Y')^3),
\end{array}
\end{equation}
\begin{equation}\label{div_x_sub}
\begin{array}{c}
\displaystyle{\rm div}_{x'}\left(\int_0^{h_{\rm max}}\tilde u'(x',y_3)\,dy_3\right)=0\hbox{  in  }\omega,\\
\noame
\displaystyle
\left(\int_0^{h_{\rm max}}\tilde u'(x',y_3)\,dy_3\right)\cdot n=0\hbox{  in  }\partial\omega\,,
\end{array}
\end{equation}
\begin{equation}
\begin{array}{c}\label{div_y_sub}
 \displaystyle{\rm div}_{y'}\hat u'=0\hbox{  in  }\omega\times Y,\\
 \noame
\displaystyle
{\rm div}_{x'}\left(\int_{Y}\hat u'(x',y)\,dy\right)=0\hbox{  in  }\omega,\  \left(\int_{Y}\hat u'(x',y)\,dy\right)\cdot n=0
\hbox{  in  }\partial\omega\,,
\end{array}
\end{equation}

\item[(ii)] (Microrotation) there exist $\tilde w\in H^1(0,h_{\rm max};L^2(\omega)^3)$, with $\tilde w=0$ on $y_3=\{0,h_{\rm max}\}$ and $\tilde w_3=0$, and $\hat w\in H^1(0,h(y'); L^2_{\#}(\omega\times Y')^3)$, with $\hat w=0$ on $y_3=\{0,h(y')\}$, such that $\int_{Y}\hat w(x',y)dy=\int_0^{h_{\rm max}}\tilde w(x',$ $y_3)\,dy_3$ with $\int_{Y}\hat w_3\,dy=0$,  and moreover 
\begin{equation}\label{conv_w_sub}
\begin{array}{c}
\displaystyle\eta_\ep^{-1}\tilde w_\ep\rightharpoonup (\tilde w,0)\hbox{  in  }H^1(0,h_{\rm max};L^2(\omega)^3),\\
\noame
\displaystyle 
\eta_\ep^{-1}\hat w_\ep\rightharpoonup \hat w\hbox{  in  }H^1(0,h(y'); L^2(\omega\times Y')^3),
\end{array}
\end{equation}
\item[(iii)] (Pressure) there exists $\tilde P\in L^2_0(\Omega)$ independent of $y_3$, such that 
\begin{eqnarray}
&\displaystyle\tilde P_\ep\to \tilde P\hbox{  in  }L^2(\Omega),\quad \hat P_\ep\to \tilde P\hbox{  in  }L^2(\omega\times \Pi).& \label{conv_P_sub}
\end{eqnarray}
\end{itemize}
\end{lemma}
\begin{proof} The proof of $(i)$ is similar to the critical case, but we have to take into account that applying the change of variables (\ref{CV}) to the divergence condition ${\rm div}_{\eta_\ep}\tilde u_\ep$, multiplying by $\eta_\ep^{-1}$ and passing to the limit, we prove that $\hat u_3$ is independent of $y_3$. Thus, the divergence condition on $y'$ given in (\ref{div_y_sub}) straightforward. For more details, we refer the reader to the proof of Lemmas 5.2-i) and 5.4-ii) in \cite{Anguiano_SG}.

The proofs of $(ii)$  and $(iii)$ are similar to the critical case, so we omit it.
\hfill$\square$
\end{proof}

Next, we give the homogenized system satisfied by $(\hat u,\hat w,\tilde P)$.
\begin{theorem}
In the case $\eta_\ep\ll \ep$, then the sequence $(\eta_\ep^{-2}\hat u_\ep, \eta_\ep^{-1}\hat w_\ep)$   converges weakly to 
$(\hat u,\hat w)$ in $H^1(0,h(y'); L^2(\omega\times Y')^3)\times H^1(0,h(y'); L^2(\omega\times Y')^3)$ and $\hat P_\ep$ converges strongly to $\tilde P$ in $L^2(\Omega)$, where  $(\hat u,\hat w, \tilde P)$ in $H^1(0,h(y'); L^2_{\#}(\omega\times Y')^3)\times H^1(0,h(y'); L^2_{\#}(\omega\times Y')^3)\times (L^2_0(\omega)\cap H^1(\omega))$ with  $\int_{Y}\hat u_3\,dy=0$, $\hat u_3$ independent of $y_3$ and $\hat w_3=0$, is the unique solution of the following homogenized system
\begin{equation}\label{hom_system_sub}
\!\!\left\{\begin{array}{rl}
\displaystyle
-\partial_{y_3}^2  \hat u'+\nabla_{y'} \hat q=2N^2 {\rm rot}_{y_3}\hat w'+f'(x')-\nabla_{x'}\tilde P(x')&\hbox{in }\omega\times Y,\\
\noame
{\rm div}_{y'}\hat u'=0&\hbox{in }\omega\times Y,\\
\noame
-R_c\partial_{y_3}^2 \hat w'+4N^2 \hat w'=2N^2 {\rm rot}_{y_3}\hat u'+g'(x')&\hbox{in }\omega\times Y,\\
\noame
\hat u'=0&\hbox{on }y_3=\{0,h(y')\},\\
\noame
\displaystyle{\rm div}_{x'}\left(\int_{Y}\hat u'(x',y)\,dy\right)=0&\hbox{in }\omega,\\
\noame
\displaystyle\left(\int_{Y}\hat u'(x',y)\,dy\right)\cdot n=0&\hbox{on }\partial\omega,\\
\noame
\hat u(x',y), \hat w(x',y), \hat q(x',y')\quad Y'-\hbox{periodic}.&
\end{array}\right.
\end{equation}
\end{theorem}
\begin{proof}From Lemma \ref{lem_asymp_sub}, conditions  (\ref{hom_system_sub})$_{2,4,5,6}$ hold. To prove that $(\hat u, \hat w, \tilde P)$ sa\-tisfies the momentum equations
given in (\ref{hom_system_sub}), we consider $\varphi\in \mathcal{D}(\omega;C_\#^\infty(Y)^3)$ with $\varphi_3$ independent of $y_3$, ${\rm div}_{y'}\varphi'=0$ in $\omega\times Y$ and ${\rm div}_{x'}\int_{Y}\varphi'\,dy=0$ in $\omega$, and we choose $\varphi_\ep=(\varphi',\varphi_3)$ in (\ref{form_var_hat_u}).

Taking into account that  ${\rm div}_{y'}\varphi'=0$ in $\omega\times Y$ and $\varphi_3$ is independent of $y_3$, we have that 
$${1\over \eta_\ep}\int_{\omega\times \Pi}\hat P_\ep\,{\rm div}_{y'}\varphi'\,dx'dy=0\quad\hbox{and}\quad{1\over \eta_\ep}\int_{\omega\times \Pi}\hat P_\ep\,\partial_{y_3}\varphi_3\,dx'dy=0\,.$$
Also, from Cauchy-Schwarz's inequality, the second estimate in (\ref{estim_u_hat}), the first estimate in (\ref{estim_w_hat}) and $\eta_\ep/\ep\to 0$, we have that 
$$\left|{1\over \varepsilon^{2}}\int_{\omega\times Y}D_{y'}\hat u_\varepsilon:D_{y'}\varphi\,dx'dy_3\right|\leq {C\over \varepsilon^2}\|D_{y'}\hat u_\varepsilon\|_{L^2(\omega\times Y)^{3\times 2}}\leq  C{\eta_\varepsilon\over \varepsilon}\to 0,$$
and 
$$\begin{array}{rl}
\displaystyle\left|{1\over \varepsilon}\int_{\omega\times Y}{\rm rot}_{y'}\hat w_{\varepsilon,3}\cdot \varphi'\,dx'dy\right| &\displaystyle=\left|{1\over \varepsilon}\int_{\omega\times Y}\hat w_{\varepsilon,3}\cdot {\rm rot}_{y'}\varphi'\,dx'dy\right|\\
\noame&\displaystyle
\leq {C\over \varepsilon}\|w_\varepsilon\|_{L^2(\omega\times Y)^3}\leq C{\eta_\varepsilon\over \varepsilon}\to 0.
\end{array}$$
Thus, passing to the limit using the convergences (\ref{conv_u_sub}), (\ref{conv_w_sub}) and (\ref{conv_P_sub}), we obtain
$$\begin{array}{l}
\displaystyle\int_{\omega\times Y}\partial_{y_3}\hat u'\cdot \partial_{y_3}\varphi'\,dx'dy-\int_{\omega\times Y}\tilde P\,{\rm div}_{x'}\varphi'\,dx'dy
\\
\noame
\displaystyle=2N^2\int_{\omega\times Y}{\rm rot}_{y_3}\hat w'\cdot \varphi'\,dx'dy+\int_{\omega\times Y}f'\cdot \varphi'\,dx'dy\,.
\end{array}$$
Since $\hat P$ does not depend on $y$ and ${\rm div}_{x'}\int_Y\varphi'\,dy=0$ in $\omega$, we have that
$$\int_{\omega\times Y}\tilde P\,{\rm div}_{x'}\varphi'\,dx'dy=\int_{\omega}\tilde P\,{\rm div}_{x'}\left(\int_Y \varphi'\,dy\right)dx'=0,$$
so we get
\begin{equation}\label{form_var_limit_sub}
\begin{array}{l}
\displaystyle\int_{\omega\times Y}\partial_{y_3}\hat u'\cdot \partial_{y_3}\varphi'\,dx'dy\\
\noame
\displaystyle=2N^2\int_{\omega\times Y}{\rm rot}_{y_3}\hat w'\cdot \varphi'\,dx'dy+\int_{\omega\times Y}f'\cdot \varphi'\,dx'dy\,.
\end{array}\end{equation}

Next, for every $\psi\in \mathcal{D}(\omega;C_\#^\infty(Y)^3)$, we choose $\psi_\ep=\eta_\ep^{-1}\psi$ in (\ref{form_var_hat_w}). Then, proceeding similarly as above and passing to the limit using previous convergences, we get
\begin{equation}\label{form_var_limit_sub_w}
\begin{array}{l}
\displaystyle
R_c\int_{\omega\times Y}\partial_{y_3}\hat w'\cdot \partial_{y_3}\psi'\,dx'dy+4N^2\int_{\omega\times Y}\hat w'\cdot \psi'\,dx'dy
\\
\noame
\displaystyle
=2N^2\int_{\omega\times Y}{\rm rot}_{y_3}\hat u'\cdot \psi'\,dx'dy+\int_{\omega\times Y}g'\cdot \psi'\,dx'dy.
\end{array}
\end{equation}
Finally, we can prove $\hat w_3=0$. For this, we take  in (\ref{form_var_hat_w}) the test function $\psi_\ep=(0,\eta_\ep^{-1}\psi_3)$, and passing to the limit as above, we get 
$$\begin{array}{l}
\displaystyle
R_c\int_{\omega\times Y}\partial_{y_3}\hat w_3:\partial_{y_3}\psi_3\,dx'dy+4N^2\int_{\omega\times Y}\hat w_3\cdot \psi_3\,dx'dy
=0,
\end{array}$$
which is equivalent to the equation $-R_c\partial_{y_3}^2\hat w_3+2N^2\hat w_3=0$. This  together with the boundary conditions $\hat w_3=0$ on $y_3=\{0,h(y')\}$ implies that $\hat w_3=0$.

By density, and reasoning as in the proof of Theorem \ref{hom_system_crit},  problem (\ref{form_var_limit_sub})-(\ref{form_var_limit_sub_w}) is equivalent to the homogenized system (\ref{hom_system_sub}) (observe that the condition ${\rm div}_{y'}\varphi'=0$ implies that $\hat q$ does not depend on $y_3$). Since $\partial_{y_3}\hat u' + 2N^2 {\rm rot}_{y_3}\hat w'+ f'\in L^2(\omega\times Y)$, it can be easily proved that $\nabla_{x'}\tilde P\in L^2(\omega)^2$ and so $\tilde P\in H^1(\omega)$ and also that system (\ref{hom_system_sub}) has a unique solution (see for example Proposition 3.3 and 3.5 in \cite{MT}).
\hfill$\square$
\end{proof}
Let us define the local problems which are useful to eliminate the variable $y$ of the previous homogenized problem and then obtain a Reynolds equation for $\tilde P$.  

We define  $\Phi$ and $\Psi$ by 
\begin{equation}\label{Phi}
\begin{array}{rl}
\displaystyle
\Phi(h(y'),N,R_c)=& \displaystyle {1\over 12}+{R_c\over 4h^2(y')(1-N^2)}\\
\noame
&\displaystyle-{1\over 4h(y')}\sqrt{{N^2 R_c\over 1-N^2}}\coth\left(Nh(y')\sqrt{{1-N^2\over R_c}}\right)\,,
\end{array}
\end{equation}
\begin{equation}\label{Psi}
\Psi(h(y'),N,R_c)={\tanh\left(Nh(y')\sqrt{{1-N^2\over R_c}}\right)\over {1-{N\over h(y')}\sqrt{{1-N^2\over R_c}}\tanh\left(Nh(y')\sqrt{{1-N^2\over R_c}}\right)}}\,,
\end{equation}
 and for every $i,k=1,2$, we consider the following local Reynolds problems 
\begin{equation}\label{local_problems_sub_2}
-{\rm div}_{y'}\left({h^3(y')\over 1-N^2}\Phi(h(y'),N,R_c)\left(\nabla_{y'}\pi^{i,k}(y')+e_i\delta_{1k}\right)\right)=0\ \hbox{  in  }Y'\,.
\end{equation}

It is known that from the positivity of function $\Phi$, problem (\ref{local_problems_sub_2}) has a unique solution for $\pi^{i,k}\in H^1_{\#}(Y')$ (see \cite{BayadaChamGam} for more details).\\

Next, we give the main result of this section. 
\begin{theorem}\label{them_main_sub}
Let $(\hat u,\hat w,\tilde P)\in L^2(\omega;H^1_\#(Y)^3)\times L^2(\omega;H^1_\#(Y)^3)\times (L_0^2(\omega)\cap H^1(\omega))$ be the unique weak solution of problem (\ref{hom_system_sub}). Then, the extensions $(\eta_\ep^{-2}\tilde u_\ep,$ $\eta_\ep^{-1}\tilde w_\ep)$ and $\tilde P_\ep$ of the solution of problem (\ref{system_2})-(\ref{bc_system_2}) converge weakly to $(\tilde u,\tilde w)$ in $H^1(0,h_{\rm max};L^2(\omega)^3)\times H^1(0,h_{\rm max};L^2(\omega)^3)$ and strongly to  $\tilde P$ in $ L^2(\Omega)$ res\-pectively, with $\tilde u_3=\tilde w_3=0$. Moreover, defining $\widetilde U(x')=\int_0^{h_{\rm max}}\tilde u(x',y_3)\,dy_3$ and $\widetilde W(x')=\int_0^{h_{\rm max}}\tilde w(x',y_3)\,dy_3$, it holds 
\begin{equation}\label{Darcy_law_u_w_sub_1}
\begin{array}{ll}
\widetilde U'(x')=K^{(1)}_0\left(f'(x')-\nabla_{x'}\tilde P(x')\right)+K^{(2)}_0 g'(x'),&\quad \widetilde U_3(x')=0\quad \hbox{in }\omega,\\
\noame
\displaystyle \widetilde W'(x')=L_0^{(2)}\,g'(x'),&\quad \widetilde W_3(x')=0\quad \hbox{in }\omega,
\end{array}
\end{equation}
where the matrices $K^{(k)}_0$, $k=1,2$,  and $L_0^{(2)}$ are matrices with coefficients 
\begin{equation}\label{def_K_sub}
\begin{array}{l}
\displaystyle \left(K^{(k)}_0\right)_{ij}={1\over 1-N^2}\int_{Y'}{h^3(y')}\Phi(h(y'),N,R_c)\left(\partial_{y_i}\pi^{j,k}(y')+\delta_{ij}\delta_{1k}\right)dy',\\
\noame
\displaystyle \left(L_0^{(2)}\right)_{ij}=-{1\over 4N^3}\sqrt{{R_c\over 1-N^2}}\left(\int_{Y'}\Psi(h(y'),N)\,dy'\right)\delta_{ij}\,,
\end{array}
\end{equation}
for $i,j=1,2$, with  $\Phi$ and $\Psi$ given by (\ref{Phi}) and (\ref{Psi}), respectively, and $\pi^{i,k}\in H^1_{\#}(Y')$, $i,k=1,2$, the unique solutions of the cell problems (\ref{local_problems_sub_2}).

Here, $\tilde P\in H^1(\omega)\cap L^2_0(\omega)$ is the unique solution of problem
\begin{equation}\label{Darcy_law_P_sub}
\left\{\begin{array}{l}
{\rm div}_{x'}\left( -A_0 \nabla_{x'}\tilde P(x')+ b_0(x')\right)=0\quad \hbox{ in }\omega,\\
\noame
\left( -A_0 \nabla_{x'}\tilde P(x')+ b_0(x')\right)\cdot n=0\quad \hbox{ on }\partial\omega,
\end{array}\right.
\end{equation}
where the flow factors are given by $A_0=K_0^{(1)}$ and $b_0(x')=K_0^{(1)}f'(x')+K_0^{(2)}g'(x')$.
\end{theorem}
\begin{proof}
We proceed as in in the proof of Theorem \ref{them_main_crit} in order to obtain (\ref{Darcy_law_u_w_sub_1}).
Thus, expressions for $\widetilde U$ and $\widetilde W$ can be obtained by defining
\begin{eqnarray}
&&\hat u'(x',y)=\sum_{i=1}^2\left[\left(\partial_{x_i}\tilde P(x')-f_i(x')\right) u^{i,1}(y)- g_i(x') \,u^{i,2}(y)\right],\nonumber\\
\noame
&&\hat w'(x',y)=\sum_{i=1}^2\left[\left(\partial_{x_i}\tilde P(x')-f_i(x')\right) w^{i,1}(y)- g_i(x')\, w^{i,2}(y)\right],\label{identifications_sub}\\
\noame
&&\hat q(x',y)=\sum_{i=1}^2\left[\left(\partial_{x_i}\tilde P(x')-f_i(x')\right) \pi^{i,1}(y')- g_i(x')\, \pi^{i,2}(y')\right],\nonumber
\end{eqnarray}
where $(u^{i,k},w^{i,k})\in H^1_\#(Y)^2\times H^1_\#(Y)^2$, $i,k=1,2$, are the unique solutions of
\begin{equation}\label{local_problems_sub_1}
\left\{\begin{array}{rl}\displaystyle
-\partial_{y_3}^2 u^{i,k}+\nabla_{y'}\pi^{i,k}-2N^2{\rm rot}_{y_3} w^{i,k}=-e_i\delta_{1k}&\hbox{ in }Y,\\
\noame
{\rm div}_{y'} u^{i,k}=0&\hbox{ in }Y,\\
\noame
-R_c\partial_{y_3}^2 w^{i,k}+4N^2 w^{i,k}-2N^2{\rm rot}_{y_3} u^{i,k}=-e_i\delta_{2k}&\hbox{ in }Y,\\
\noame
u^{i,k}=w^{i,k}=0&\hbox{ on }y_3=\{0,h(y')\},\\
\noame
u^{i,k}(y), w^{i,k}(y),\pi^{i,k}(y')\quad Y'-\hbox{periodic}\,.
\end{array}\right.
\end{equation}
Then, thanks to the identities $\int_Y \hat u'(x',y)\,dy=\int_0^{h_{\rm max}}\tilde u'(x',y_3)\,dy_3$,  $\int_Y\hat u_3\,dy=0$,   $\int_Y\hat w'(x',y)\,dy=\int_0^{h_{\rm max}}\tilde w'(x',y_3)\,dy_3$ and 
$\hat w_3=0$ given in Lemma \ref{lem_asymp_sub}, it holds 
\begin{equation}\label{Darcy_law_u_w_sub_proof}
\begin{array}{l}
\widetilde U'(x')=\displaystyle \int_Y\hat u'(x',y)\,dy=-K^{(1)}_0\left(\nabla_{x'}\tilde P(x')-f'(x')\right)+K^{(2)}_0 g'(x'),\\
\noame
\displaystyle \widetilde U_3(x')=\displaystyle \int_{Y}\hat u_3(x',y')\,dy=0\quad \hbox{ in }\omega,\\
\noame
\displaystyle \widetilde W'(x')=\displaystyle \int_Y\hat w'(x',y)\,dy=-L^{(1)}_0\left(\nabla_{x'}\tilde P(x')-f'(x')\right)+L^{(2)}_0 g'(x'),
\\
\noame
\displaystyle \widetilde W_3(x')=\displaystyle \int_Y\hat w_3(x',y)\,dy=0\quad \hbox{ in }\omega,
\end{array}
\end{equation}
where $K^{(k)}_0$, $L^{(k)}_0$, $k=1,2$, are matrices defined by their coefficients
\begin{equation}\label{K_0_proof} \left(K^{(k)}_0\right)_{ij}=-\int_{Y}u^{i,k}_j(y)\,dy,\quad \left(L^{(k)}_0\right)_{ij}=-\int_{Y}w^{i,k}_j(y)\,dy\,,\quad i,j=1,2\,.
\end{equation}
Then, by the divergence condition in the variable $x'$ given in (\ref{hom_system_sub}), we get the generalized Reynolds equation (\ref{Darcy_law_P_sub}).

However, we observe that (\ref{local_problems_sub_1}) can be viewed as a system of ordinary diffe\-rential equations with constant coefficients, with respect to the variable $y_3$ and unkowns functions $y_3\mapsto u^{i,k}_1 (y',y_3),w^{i,k}_2 (y',y_3),  u^{i,k}_2 (y',y_3), w^{i,k}_1 (y',y_3)$, where $y'$ is a parameter, $y'\in Y'$. Thus, we can give explicit expressions for $u^{i,k}$ and $w^{i,k}$. 

The procedure to obtain a solution to the previous system  is given in the Appendix (see also in \cite{BayadaChamGam,BayadaLuc}). Thus, considering $\bar u'=u^{i,k}$, $\bar w'=w^{i,k}$, $\bar f'=-e_i\delta_{ik}$ and  $\bar g'=-e_i\delta_{2k}$ in (\ref{expressions_final_app})-(\ref{expressions_final_app123}), we obtain that $u^{i,k},w^{i,k}$ are given in terms of $\pi^{i,k}$ by the expressions
\begin{equation}\label{expression_cell_u_w}\begin{array}{l}
u^{i,k}(y)=
 {1\over 2(1-N^2)}\Big[y_3^2-h(y')y_3 \\
 \noame   +{h(y')N^2\over k}\left(\sinh({ky_3})-(\cosh(ky_3)-1)\coth{\left({k h(y')\over 2}\right)}\right)\Big]\left(\nabla_{y'}\pi^{i,k}(y')+e_i\delta_{1k}\right)\\
\noame 
 +{h(y')\over N^2}\left[\left({2N^2\over k}\sinh(ky_3)-2y_3\right)A+{2N^2\over k}(\cosh(ky_3)-1)B-y_3\right]\left(e_i\delta_{2k}\right)^{\perp}\,,
 \end{array}\end{equation}
$$\begin{array}{l}
w^{i,k}(y)= {1\over 4(1-N^2)}\Big[2y_3\\
\noame
+h(y')\left(\cosh(ky_3)-1-\sinh{(ky_3)}\coth\left({k h(y')\over 2}\right)\right)\Big]\left(\nabla_{y'}\pi^{i,k}(y')+e_i\delta_{1k}\right)^\perp\\
\noame
 -{h(y')\over 2N^2}\Big[\cosh(ky_3)A+\sinh(ky_3)B\Big]e_i\delta_{2k}\,,\\
\end{array}
$$
where $k=\sqrt{{4N^2(1-N^2)\over R_c}}$ and $A$, $B$ are given by
\begin{equation}\label{def_AB}
\begin{array}{l}    A(y')={\sinh(kh(y'))\over -2h(y')\sinh(kh(y'))+{4N^2\over k}(\cosh(kh(y'))-1)},\\
\noame
 B(y')={-(\cosh(kh(y')-1)\over -2h(y')\sinh(kh(y'))+{4N^2\over k}(\cosh(kh(y'))-1)}\,.
 \end{array}
 \end{equation}
Taking into account that from (\ref{int_u_w_app}) it holds 
\begin{equation}\label{property_int_uik}
\begin{array}{l}\displaystyle
\int_0^{h(y')}u^{i,k}(y',y_3)\,dy_3=-{h^3(y')\over 1-N^2}\Phi(h(y'),N,R_c)\left(\nabla_{y'}\pi^{i,k}+e_i\delta_{1k}\right)\,,\\
\noame
\displaystyle\int_0^{h(y')}w^{i,k}(y',y_3)\,dy_3=-{1\over 4N^3}\sqrt{{R_c\over 1-N^2}}\Psi(h(y'),N,R_c)e_i\delta_{2k}\,,
\end{array}
\end{equation}
with $\Phi$ and $\Psi$ given by (\ref{Phi}) and (\ref{Psi}), respectively, we get that $\pi^{i,k}$ satisfies the generalized Reynolds cell problem (\ref{local_problems_sub_2}).  Using the expressions of $u^{i,k}$ and $w^{i,k}$ together with (\ref{Darcy_law_u_w_sub_proof}), (\ref{K_0_proof}) and (\ref{property_int_uik}), we easily get (\ref{Darcy_law_u_w_sub_1}). Observe that, from the second equation in (\ref{property_int_uik}) with $k=2$, we have $L^{(1)}_0=0$, which ends the proof. 
\hfill$\square$\end{proof}

\section{High-frequency roughness regime ($\lambda=+\infty$)}\label{HFRR}
It corresponds to the case when the wavelength of the roughness is much smaller than the film thickness, i.e. $\eta_\ep\gg \ep$ which is equivalent to $\lambda=+\infty$.

Next, we give some compactness results about the behavior of the extended sequence $(\tilde u_\ep,\tilde w_\ep,\tilde P_\ep)$ and the unfolding functions $(\hat u_\ep,\hat w_\ep, \hat P_\ep)$ satisfying the {\it a priori} estimates given in Lemmas \ref{lemma_estimates} and \ref{lemma_est_P}, and Lemma \ref{estimates_hat}, respectively.

\begin{lemma}\label{lem_asymp_sup}
For a subsequence of $\ep$ still denoted by $\ep$, we have that
\begin{itemize}
\item[(i)] (Velocity) there exists  $\tilde u\in H^1(0,h_{\rm max};L^2(\omega)^3)$, with $\tilde u=0$ on $y_3=\{0,h_{\rm max}\}$ and $\tilde u_3=0$, such that 
\begin{equation}\label{conv_u_sup}
\begin{array}{c}
\displaystyle\eta_\ep^{-2}\tilde u_\ep\rightharpoonup (\tilde u',0)\hbox{  in  }H^1(0,h_{\rm max};L^2(\omega)^3),\\
\noame
\eta_\ep^{-2}\tilde u_\ep\rightharpoonup 0\hbox{  in  }H^1(h_{\rm min},h_{\rm max};L^2(\omega)^3),
\end{array}
\end{equation}
\begin{equation}\label{conv_u_sup_1}
\displaystyle\eta_\ep^{-2}\hat u_\ep\rightharpoonup (\tilde u',0)\hbox{  in  }H^1(0,h_{\rm min};L^2(\omega)^3),
\end{equation}
\begin{equation}\label{div_x_sup}
\begin{array}{c}
\displaystyle
{\rm div}_{x'}\left(\int_0^{h_{\rm min}}\tilde u'(x',y_3)\,dy_3\right)=0\hbox{  in  }\omega,\\
\noame
\displaystyle
\left(\int_0^{h_{\rm min}}\tilde u'(x',y_3)\,dy_3\right)\cdot n=0\hbox{  in  }\partial\omega\,,
\end{array}\end{equation}

\item[(ii)] (Microrotation) there exists  $\tilde w\in H^1(0,h_{\rm max};L^2(\omega)^3)$, with $\tilde w=0$ on $y_3=\{0,h_{\rm max}\}$ and $\tilde w_3=0$, such that 
\begin{equation}\label{conv_w_sup}
\begin{array}{c}
\displaystyle\eta_\ep^{-1}\tilde w_\ep\rightharpoonup (\tilde w',0)\hbox{  in  }H^1(0,h_{\rm max};L^2(\omega)^3),\\
\noame
\displaystyle
\eta_\ep^{-1}\tilde w_\ep\rightharpoonup 0\hbox{  in  }H^1(h_{\rm min},h_{\rm max};L^2(\omega)^3),
\end{array}\end{equation}
\begin{equation}\label{conv_w_sup_1}
\displaystyle\eta_\ep^{-1}\hat w_\ep\rightharpoonup (\tilde w',0)\hbox{  in  }H^1(0,h_{\rm min};L^2(\omega)^3),
\end{equation}
\item[(iii)] (Pressure) there exists a function $\tilde P\in L^2_0(\Omega)$ independent of $y_3$, such that
\begin{eqnarray}
&\displaystyle\tilde P_\ep\to \tilde P\hbox{  in  }L^2(\Omega),\quad \hat P_\ep\to \tilde P\hbox{  in  }L^2(\omega\times \Pi)\,.& \label{conv_P_sup}
\end{eqnarray}
\end{itemize}
\end{lemma}
\begin{proof} We start proving $(i)$. We will only give some remarks and for more details, we refer to the reader to  Lemmas 5.2-ii) and 5.4-ii) in \cite{Anguiano_SG}. As pre\-vious cases, we can prove that there exists $\tilde u\in H^1(0,h_{\rm max};L^2(\omega)^3)$ such that $\eta_\ep^{-2}\tilde u_\ep$ converges weakly to $\tilde u$ in $H^1(0,h_{\rm max};L^2(\omega)^3)$. On the other hand, from estimate (\ref{estim_tilde_u_ep+}),  $\ep/\eta_\ep\to 0$ and taking into account that $\eta_\ep^{-2}=({\ep\over \eta_\ep})^2\ep^{-2}$, then second convergence in (\ref{conv_u_sup}) holds and so $\tilde u=0$ in $\Omega^+$. Then, reasoning as previous cases, we can prove that $\tilde u_3=0$, $\tilde v'=0$ on $y_3=\{0,h_{\rm min}\}$ and also, the divergence condition (\ref{div_x_sup}). 

From estimates (\ref{estim_u_hat}), we deduce that  there exists $\hat u\in H^1(0,h(y');L^2(\omega\times Y')^3)$ such that 
\begin{equation}\label{conv_u_hat_proof}
\hat u_\ep\rightharpoonup \hat u\hbox{   in   } H^1(0,h(y');L^2(\omega\times Y')^3).
\end{equation}  Since $\ep^{-1}\eta_\ep^{-1}D_y\hat u_\ep$ is bounded in $L^2(\omega\times Y)^3$, we observe that $\eta_\ep^{-2}D_y\hat u_\ep$ is also bounded, and tends to zero. This together with (\ref{conv_u_hat_proof}) implies $\eta_\ep^{-2}D_{y'}\hat u_\ep$ converges weakly to zero in $H^1(0,h(y');L^2(\omega\times Y')^{3\times 2})$, and so $\hat u$ does not depend on $y'$.  

Proceeding as previous cases, but taking $\varphi\in C^1_c(\Omega^+)$, we can prove that
$$\int_{\omega\times \Pi^+}\hat u(x',y)\varphi(x',y_3)\,dy=\int_{\Omega^+}\tilde u(x',y_3)\varphi(x',y_3)\,dx'dy_3,$$
and taking into account that $\tilde u=0$ on $\Omega^+$, we deduce that $\hat u=0$ in $\omega\times \Pi^+$. Then, we can prove that $\int_{\omega\times \Pi^-}\hat u(x',y)\varphi(x',y_3)\,dy=\int_{\Omega^-}\tilde u(x',y_3)\varphi(x',y_3)$ $dx'dy_3$ holds and, since $\hat u$ does not depend on $y'$, we have that $\hat u=(\tilde u',0)$.

For the proof of $(ii)$ for microrotation, we can proceed as for the velo\-city.  By considering estimate (\ref{estim_tilde_w_ep}), we prove the existence of the weak limit $\tilde w\in H^1(0,h(y');L^2(\omega\times Y')^3) $ of the sequence $\eta_\ep^{-1}\tilde w_\ep$, and taking into account estimate (\ref{estim_tilde_w_ep+}), $\ep/\eta_\ep\to 0$ and  that $\eta_\ep^{-1}=(\eta_\ep^2\ep^3)({\ep\over \eta_\ep})^3$, we prove the second convergence in (\ref{conv_w_sup}). Moreover, as in the case of the velocity, it can be proved that $\tilde w=0$ on $y_3=\{0,h_{\rm min}\}$.  To prove that $\tilde w_3=0$, we argue as in the critical case, by taking a test function $\psi_\ep=(0,0,\eta_\ep^{-1}\psi_3)$  in (\ref{form_var_general_1}),  passing to the limit and considering the previous boundary conditions. For the proof of (\ref{conv_w_sup_1}), we proceed as the case of the velocity by taking into account estimates (\ref{estim_w_hat}).

Finally, to prove $(iii)$, we proceed as in the critical case. First we prove weak convergence of the extended pressure $\tilde P_\ep$ to a  function $\tilde P$ in $L^2_0(\Omega)$ and next,  we prove that $\tilde P$ independent of $y_3$. Finally,    we prove strong convergence of the pressure, but in this case  we have to take into account te behavior of $\tilde u_\ep$ and $\tilde w_\ep$ on the oscillating part.  Thus, we consider $\sigma_\varepsilon\in H^1_0(\Omega)^3$ such that $\sigma_\ep\rightharpoonup \sigma$ in $H^1_0(\Omega)^3$.  Denoting $\tilde\sigma_\varepsilon=(\sigma_\varepsilon',\varepsilon \sigma_{\varepsilon,3})$ and $\tilde\sigma=(\sigma',0)$, we have 
\begin{equation}\label{convsigmappp}
\tilde\sigma_\varepsilon\rightharpoonup \tilde \sigma\quad\hbox{in }H^1_0(\Omega)^3.
\end{equation}
Then,
$$\begin{array}{l}
\displaystyle\left|<\nabla_{x',y_3}\tilde P_\varepsilon,\sigma_\varepsilon>_{\Omega^+}-<\nabla_{x'}\tilde P,\tilde \sigma>_{\Omega^+}\right|
\\
\noame
\leq
\displaystyle\left|<\nabla_{x',y_3}\tilde P_\varepsilon-\nabla_{x'}\tilde P,\tilde \sigma>_{\Omega^+}\right|
+\left|<\nabla_{x',y_3}\tilde P_\varepsilon,\sigma_\varepsilon-\tilde\sigma>_{\Omega^+}\right|\,.
\end{array}
$$
On the one hand, using the weak convergence of the pressure, we have 
$$\left|<\nabla_{x',y_3} \tilde P_\ep-\nabla_{x'}\tilde P,\tilde\sigma>_{\Omega^+}\right|=\left|\int_{\Omega^+}\left(\tilde P_\ep-\tilde P\right)\,{\rm div}_{x'}\sigma'\,dx\right|\to 0,\quad \hbox{as }\ep\to 0\,.$$
On the other hand, proceeding as in Lemma \ref{lemma_est_P}, we have that
$$\begin{array}{l}
\left|<\nabla_{x',y_3}\tilde P_\ep,\sigma_\ep-\tilde \sigma>_{\Omega^+}\right|=\left|<\nabla_{\eta_\varepsilon}\tilde P_\ep,\tilde \sigma_\ep-\tilde \sigma)>_{ \Omega^+}\right|\\
\noame
\quad \displaystyle
\leq C\left(\|\tilde\sigma_\varepsilon-\tilde\sigma\|_{L^2(\Omega)^3}+ \varepsilon\|D_{x',y_3}(\sigma_\varepsilon-\tilde \sigma)\|_{L^2(\Omega)^{3\times 3}}\right),
\end{array}$$
which tends to zero because of the convergence of the sequence $\sigma_\ep$ and the Rellich theorem. 

Then, reasoning similarly as above by considering in $\Omega^-$, we deduce that 
$$\left|<\nabla_{x',y_3}\tilde P_\ep,\sigma_\ep>_{\Omega^-}-<\nabla_{x'}\tilde P,\sigma>_{\Omega^-}\right|\to 0\,,$$
which together with previous convergence, implies the convergence of $\nabla_{x',y_3}\tilde P_\ep$ to $\nabla_{x'}\tilde P$ strongly in $H^{-1}(\Omega)^3$.  This together with the  Ne${\breve{\rm c}}$as inequality (\ref{Necas_inq}) implies  the first  convergence in (\ref{conv_P_sup}). Finally, we remark that the strong convergence of sequence $\hat P_\ep$ to $\tilde P$ is a consequence of the strong convergence of $\tilde P_\ep$ to $\tilde P$ (see  \cite[Proposition 2.9]{Ciora2}).
\hfill$\square$
\end{proof}

As seen in the previous compactness result, the microstructure of $\widetilde\Omega_\ep$ will not be involved in the homogenized system and thus, we will obtain a Reynolds equation satisfied by $\tilde P$ in the non-oscillating part of the domain, that is $\Omega^-$.
\begin{theorem}\label{thm_general_sup}
In the case $\eta_\ep\gg \ep$, then the extensions $(\eta_\ep^{-2}\tilde u_\ep, \eta_\ep^{-1}\tilde w_\ep)$ and  $\tilde P_\ep$ of the solution of problem (\ref{system_2})-(\ref{bc_system_2}) converge weakly to $H^1(0, h_{\rm min};L^2(\omega)^3)\times H^1(0, h_{\rm min};L^2(\omega)^3)$ and strongly to  $\tilde P$ in $L^2(\Omega)$ respectively, with $\tilde u_3=\tilde w_3=0$, where  $\tilde u'$ and $\tilde w'$ are given by the following expressions  in terms of the pressure $\tilde P$ in $\Omega^-$,
\begin{eqnarray}
&&\tilde u'(x',y_3)= \left[{y_3^2\over 2(1-N^2)}+{1\over 4(1-N^2)}\left({2N^2\over k}\sinh(ky_3)-2y_3\right)\right.\nonumber\\
\noame
&&\left.-{h_{\rm min}\over 2(1-N^2)}{N^2\over k}(\cosh(kh_{\rm min})-1)\coth\left({kh_{\rm min}\over 2}\right)\right]\left(\nabla_{x'}\tilde P(x')-f'(x')\right)\nonumber
\\
\noame
&&\bar w'(x',y_3)=\Big[{y_3\over 2(1-N^2)} +{h_{\rm min}\over 4(1-N^2)}\Big(\cosh(ky_3)-1\label{expressions_final_app_sup}\\
\noame
&&-\coth\left({kh_{\rm min}\over 2}\right)\sinh(ky_3)\Big)\Big]\left(\nabla_{x'}\tilde P(x')-f'(x')\right)^\perp,\nonumber
\end{eqnarray}
with $k=\sqrt{{4N^2(1-N^2)\over R_c}}$.
Moreover, defining $\widetilde U(x')=\int_0^{h_{\rm min}}\tilde u(x',y_3)\,dy_3$ and $\widetilde W(x')=\int_0^{h_{\rm min}}\tilde w(x',y_3)\,dy_3$, it holds 
\begin{equation}\label{Darcy_law_u_w_sup}
\begin{array}{ll}
 \widetilde U'(x')={h_{\rm min}\over 1-N^2}\Phi(h_{\rm min},N,R_c)\left(f'(x')-\nabla_{x'} \tilde P(x')\right),& 
\widetilde U_3(x')=0\  \hbox{in }\omega,\\
\noame
\displaystyle \widetilde W'(x')=0,&  \widetilde W_3(x')=0\ \hbox{in }\omega,
\end{array}
\end{equation}
where $\Phi$ is given by (\ref{Phi}), and $\tilde P\in H^1(\omega)\times L^2_0(\omega)$ is the unique solution of the Reynolds problem
\begin{equation}\label{Reynolds_sup_omega_minus}
\left\{\begin{array}{rl}
\displaystyle {\rm div}_{x'}\left(-A_\infty\nabla_{x'} \tilde P(x')+ b_\infty(x')\right)=0& \hbox{in }\omega\,,\\
\noame
\displaystyle \left(-A_\infty\nabla_{x'} \tilde P(x')+ b_\infty(x')\right)\cdot n=0& \hbox{on }\partial\omega\,.
\end{array}\right.
\end{equation}
Here, the flow factors are  given by $A_\infty={h_{\rm min}\over 1-N^2}\Phi(h_{\rm min},N,R_c)$ and $b_\infty(x')={h_{\rm min}\over 1-N^2}\Phi(h_{\rm min},N,R_c)f'(x')$.
\end{theorem}
\begin{proof}From Lemma \ref{lem_asymp_sup}, we observe that at main order, the microstructure does not appear because the high oscillation of the boundary. Thus, we choose in the first equation of the variational formulation (\ref{form_var_general_1}), extended to $\Omega$, the following test function $\varphi_\ep(x',y_3)=(\varphi'(x',y_3),0)\in \mathcal{D}(\Omega^-)^3$ satisfying the divergence condition ${\rm div}_{x'}\int_0^{h_{\rm min}}\varphi'(x',y_3)\,dy_3=0$ in $\omega$.
Passing to the limit by using convergences (\ref{conv_u_sup}), (\ref{conv_w_sup}) and (\ref{conv_P_sup}), we get
$$\begin{array}{l}\displaystyle\int_{\Omega^-}\partial_{y_3}\tilde u' \cdot \partial_{y_3}\varphi'\,dx'dy_3+ \int_{\Omega^-}\tilde P\,{\rm div}_{x'}\varphi'\,dx'dy_3
\\
\noame
\displaystyle= 2N^2\int_{\Omega^-}{\rm rot}_{y_3}\tilde w' \cdot \varphi'\,dx'dy_3 + \int_{\Omega^-}f'(x')\cdot \varphi'.
\end{array}$$
Since $\tilde P$ does not depend on $y_3$ and the divergence condition on the variable $x'$ satisfied by $\varphi'$, we have that 
$$\int_{\Omega^-}\tilde P\,{\rm div}_{x'}\varphi'\,dx'dy_3=\int_{\omega}\tilde P\,{\rm div}_{x'}\left(\int_0^{h_{\rm min}}\varphi'\,dy_3\right)dx'=0,$$
and so
$$\int_{\Omega^-}\partial_{y_3}\tilde u' \cdot \partial_{y_3}\varphi'\,dx'dy_3= 2N^2\int_{\Omega^-}{\rm rot}_{y_3}\tilde w' \cdot \varphi'\,dx'dy_3 + \int_{\Omega^-}f'(x')\cdot \varphi'.$$
Next, we choose in the second equation of of the variational formulation (\ref{form_var_general_1}), extended to $\Omega$, the following test function $\psi_\ep(x',y_3)=(\eta_\ep^{-1}\psi'(x',y_3),0)\in \mathcal{D}(\Omega^-)^3$ and taking into account that $\ep/\eta_\ep\to 0$, we pass to the limit and we get
$$R_c\int_{\Omega^-}\partial_{y_3}\tilde w' \cdot \partial_{y_3}\psi'\,dx'dy_3+4N^2\!\!\int_{\Omega^-}\tilde w'\cdot \psi'\,dx'dy_3= 2N^2\!\!\int_{\Omega^-}{\rm rot}_{y_3}\tilde u' \cdot \varphi'\,dx'dy_3.$$
By density arguments, previous variational formulations are equivalent to the following simplified micropolar system 
\begin{equation}\label{hom_system_sup}
\left\{\begin{array}{rl}
\displaystyle
-\partial_{y_3}^2  \tilde u'+\nabla_{x'}\tilde P(x')=2N^2 {\rm rot}_{y_3}\tilde w'+f'(x')&\hbox{ in }\Omega^-,\\
\noame
{\rm div}_{x'}\tilde u'=0&\hbox{ in }\Omega^-,\\
\noame
-R_c\partial_{y_3}^2 \tilde w'+4N^2 \tilde w'=2N^2 {\rm rot}_{y_3}\tilde u'&\hbox{ in }\Omega^-,\\
\noame
\tilde u'=0&\hbox{ on }y_3=\{0,h_{\rm hmin}\},\\
\noame
\displaystyle{\rm div}_{x'}\left(\int_0^{h_{\rm min}}\tilde u'(x',y_3)\,dy_3\right)=0&\hbox{ in }\omega,\\
\noame
\displaystyle\left(\int_0^{h_{\rm min}}\tilde u'(x',y_3)\,dy_3\right)\cdot n=0&\hbox{ on }\partial\omega\,.\end{array}\right.
\end{equation}
The solution of this system is obtained in the Appendix.  By choosing $\bar u'=\tilde u'$, $\bar w'=\tilde w'$, $\bar P=\tilde P$, $\bar f'= f'$, $\bar g'=0$ and $h(y')=h_{\rm min}$, we get expressions (\ref{expressions_final_app_sup}). By taking into account (\ref{int_u_w_app}), we get (\ref{Darcy_law_u_w_sup}), which together with the divergence condition in the variable $x'$ given in (\ref{hom_system_sup}) gives the Reynolds equation for $\tilde P$ given by (\ref{Reynolds_sup_omega_minus}).  Since $\partial_{y_3}\tilde u'\in L^2(\Omega^-)^2$, ${\rm rot}_{y_3}\tilde w'\in L^2(\Omega^-)^2$ and $f'\in L^2(\omega)$,  it can be easily proved that $\nabla_{x'}\tilde P\in L^2(\omega)^2$, and so $\tilde P\in H^1(\omega)$ and also that system (\ref{hom_system_sup}) has a unique solution (see for example Proposition 3.3 and 2.5 in \cite{MT}).
\end{proof}

\section{Conclusions}
Whereas the  multiscale analysis is well established in the lubrication field to derive a generalized  equation of the classical Reynolds equation when the boundary of the domain have small periodic oscillations, this is not the case for micropolar flows. By using dimension reduction and homogenization techniques, we studied the asymptotic behavior of the velocity, the microrotation and the pressure for a micropolar flow in a thin domain with rapidly oscillating thickness depending on two small parameters, $\eta_\ep$ and $\ep$, where $\eta_\ep$ represents the thickness of the domain and $\ep$ the wavelength of the roughness. We provide a general classification of the roughness regime for micropolar flows depending on the value $\lambda$ of the limit of  $\eta_\ep/\ep$ when $\ep$ tends to zero, which agrees with the classification of the roughness regimes for Newtonian and non-Newtonian (power law) fluids:  Stokes roughness regime ($0<\lambda<+\infty$), Reynolds roughness regime ($\lambda=0$) and high-frequency regime ($\lambda=+\infty$).
Thus, we derive three different problems,  (\ref{Darcy_law_u_w_crit})-(\ref{Darcy_law_P_crit}), (\ref{def_K_sub})-(\ref{Darcy_law_P_sub}), and (\ref{Darcy_law_u_w_sup})-(\ref{Reynolds_sup_omega_minus}), which  are written, for $0\leq \lambda\leq +\infty$, as a Reynolds equation of the form 
\begin{equation}\label{conclusion1}
\left\{\begin{array}{l}
\widetilde U'(x')=K_\lambda^{(1)}\left(f'(x')-\nabla_{x'}\tilde P(x')\right)+ K_\lambda^{(2)}g'(x'),\quad \widetilde U_3=0\hbox{ in }\omega,\\
\noame
\widetilde W'(x')=L_\lambda^{(1)}\left(f'(x')-\nabla_{x'}\tilde P(x')\right)+ L_\lambda^{(2)}g'(x'),\quad \widetilde W_3=0\hbox{ in }\omega,\\
\noame
{\rm div}_{x'}\widetilde U'(x')=0\hbox{ in }\omega,\\
\noame
 \widetilde U'(x')\cdot n=0\hbox{ on }\partial\omega\,.
\end{array}\right.
\end{equation}
 The average   velocity  $\tilde U(x')=(\widetilde U'(x'), \widetilde U_3(x'))$ and the averaged microrotation  $\tilde W(x')=(\widetilde W'(x'), \widetilde W_3(x'))$ are respectively defined by the the functions $\widetilde U(x')=\int_0^{h{\rm max}}\tilde u(x',y_3)\,dy_3$ and $\widetilde W(x')=\int_0^{h{\rm max}}\tilde w(x',y_3)\,dy_3$. We remark that in all three cases, the vertical components $\widetilde U_3$ and $\widetilde W_3$ are equal to zero.

 We observe that in (\ref{conclusion1}), $K_\lambda^{(k)}, L^{(k)}_\lambda$, $k=1,2$, $0\leq \lambda\leq +\infty$, are computed as follows:
\begin{itemize}
\item[--] In the Stokes roughness regime, $0<\lambda<+\infty$, then $K_\lambda^{(k)}, L^{(k)}_\lambda$, $k=1,2$, are calculated by solving 3D local micropolar Stokes-like problems depending on the parameter $\lambda$. We remark that the interaction between the velocity and the microrotation fields is preserved. 
\item[--] In the Reynolds roughness regime, $\lambda=0$, then $L^{(1)}_0=0$, and $K_0^{(k)}, L^{(2)}_0$, $k=1,2$, are calculated by solving 2D micropolar Reynolds-like local problems, which represents a considerable simplification.  In this case, the interaction between the velocity and the microrotation fields is also preserved. 
\item[--] In the high-frequency roughness regime, $\lambda=+\infty$, then the velocity and microrotation vanish in the oscillating zone due to the high oscillating boundary, and so we derive the classical micropolar Reynolds equation in the non-oscillating zone, where the thickness is fixed and is given by the minimum of $h$. We observe the interaction between velocity and microrotation fields is not preserved in the limit problem because only $K_\infty^{(1)}\neq 0$.
\end{itemize}
To conclude, we believe that the presented result could be instrumental for understanding the effects of the rough boundary and fluid microstructure
on the lubrication process. In view of that, more efficient numerical algorithms could be developed improving, hopefully,
the known engineering practice.

\section*{Appendix: computation of the coefficients of the micropolar Reynolds equation}
In this Appendix we describe how to obtain the coefficient of the Reynolds equation 
\begin{equation}\label{Reynolds_app}
{\rm div}_{z'}\left(-{h^3(z')\over1-N^2}\Phi(h(z'),N,R_c)\nabla_{z'}\bar p(z') + b(z')\right)=0\quad\hbox{in }\omega,
\end{equation}
where $b(x')= {h^3(z')\over1-N^2}\Phi(h(z'),N,R_c)f'(z')$ and $\Phi$ defined by (\ref{Phi}), from the micropolar system posed in $\Omega=\{(z',z_3)\in\mathbb{R}^2\times \mathbb{R}\,:\,z'\in\omega,\ 0<z_3<h(z')\}$, given by
\begin{equation}\label{problems_appendix_1}
\left\{\begin{array}{rl}\displaystyle
-\partial_{z_3}^2 \bar u_1+\partial_{z_1}\bar p(z')+2N^2\partial_{z_3} \bar w_2=\bar f_1(z')&\ \hbox{ in }\Omega,\\
\noame
-R_c\partial_{z_3}^2 \bar w_2+4N^2 \bar w_2-2N^2\partial_{z_3}\bar  u_1=\bar g_2(z')&\ \hbox{ in }\Omega,
\end{array}\right.
\end{equation}
\begin{equation}\label{problems_appendix_2}
\left\{\begin{array}{r}\displaystyle
-\partial_{z_3}^2 \bar u_2+\partial_{z_2}\bar p(z')-2N^2\partial_{z_3} \bar w_1=\bar f_2(z')\ \hbox{ in }\Omega,\\
\noame
-R_c\partial_{z_3}^2 \bar w_1+4N^2 \bar w_1+2N^2\partial_{z_3}\bar  u_2=\bar g_1(z')\ \hbox{ in }\Omega,
\end{array}\right.
\end{equation}
together with 
\begin{equation}\label{problems_appendix_3}\partial_{z_1}\left(\int_0^{h(z')}\bar u_1(z',z_3)\,dz_3\right)+ \partial_{z_2}\left(\int_0^{h(z')}\bar u_2(z',z_3)\,dz_3\right)=0\quad\hbox{ in }\omega\,,
\end{equation} 
and boundary conditions 
\begin{equation}\label{bc_app}
\bar u'(z',0)=\bar u'(z',h(z'))=\bar w'(z',0)=\bar w'(z',h(z'))=0\,.
\end{equation}

We note that $(\bar u_1, \bar w_2)$, with external forces $(\bar f',\bar g')$,  and $(\bar u_2, -\bar w_1)$, with external forces $(\bar f',-\bar g')$,  satisfy the same equations and boundary conditions. So we only describe the computation of $(\bar u_1, \bar w_2)$.\\

First, from the first equation of (\ref{problems_appendix_1}) we have
\begin{equation}\label{dem_1}
\partial_{z_3}\bar u_1(z)=\left(\partial_{z_1}\bar p(z')-\bar f_1(z')\right)z_3+2N^2\bar w_2(z) + C(z').
\end{equation}
Putting this into the second equation of (\ref{problems_appendix_1}), we have
\begin{equation}\label{dem_1}
\begin{array}{l}
\partial_{z_3}^2\bar w_2(z)-{4N^2\over R_c}(1-N^2)\bar w_2(z)\\
\noame
=-{2N^2\over R_c}\left(\partial_{z_1}\bar p(z')-\bar f_1(z')\right)z_3 -{1\over R_c} \bar g_2(z')+ {2N^2\over R_C} C(z')\,.
\end{array}
\end{equation}
The solution is 
\begin{equation}\label{w_dem}
\begin{array}{rl}
\bar w_2(z)=&A(z')\cosh(kz_3)+B(z')\sinh(kz_3)\\
\noame
&+{1\over 2(1-N^2)}(\partial_{z_1}\bar p(z')-\bar f_1(z'))z_3\\
\noame
&+{1\over 2(1-N^2)}C(z')
+{1\over 4N^2(1-N^2)}\bar g_2(z')\,,
\end{array}
\end{equation}
where $k=\sqrt{{4N^2(1-N^2)\over R_c}}$ and $A$ and $B$ are unknowns functions.

Putting this solution into equation (\ref{dem_1}), we can write $\bar u_1$ as follows
\begin{equation}\label{u_dem}
\begin{array}{rl}
\bar u_1(z)=& {z_3^2\over 2(1-N^2)}(\partial_{z_1}\bar p(z')-f_1(z'))\\
\noame
&+{2N^2\over k}(A(z')\sinh(kz_3)+B(z')\cosh(kz_3))\\
\noame & +{z_3\over 1-N^2}C(z')+{z_3\over 2(1-N^2)}\bar g_2(z')+D(z')\,.
\end{array}
\end{equation}
We rewrite $C, D$, as a function of $A$ and $B$, using the boundary conditions.
So, for $\bar u_1(z',0)=\bar w_2(z',0)=0$, we respectively get
$$D(z')=-{2N^2\over k}B(z'),\quad C(z')=2(1-N^2)\left(-A(z')-{1\over 4N^2(1-N^2)}\bar g_2(z')\right),$$
and so
\begin{equation}\label{sol_1}
\begin{array}{rl}
\bar u_1(z)=&{z_3^2\over 2(1-N^2)}(\partial_{z_1}\bar p(z')-\bar f_1(z'))+\left({2N^2\over k}\sinh(kz_3)-2z_3\right)A(z')\\
\noame
& +{2N^2\over k}(\cosh(kz_3)-1)B(z')-{z_3\over 2N^2}\bar g_2(z')\,,\\
\noame
\bar w_2(z)=&{z_3\over 2(1-N^2)}(\partial_{z_1}\bar p(z')-\bar f_1(z'))\\
\noame
&+(\cosh(kz_3)-1)A(z')+\sinh(k z_3)B(z')\,.
\end{array}
\end{equation}
Using the boundary conditions $\bar u_1(z',h(z'))=\bar w_2(z',h(z'))=0$ we get the following system
$$Q\left(\begin{array}{c}
A\\
\noame
B
\end{array}\right)=-{h(z')\over 2(1-N^2)}(\partial_{z_1}\bar p(z')-\bar f_1(z'))\left(\begin{array}{c}
h(z')\\
\noame
1
\end{array}\right)+
\bar g_2(z'){h(z')\over 2N^2} 
\left(\begin{array}{c}
1\\
\noame
0
\end{array}\right)\,,$$
where $Q$ is the matrix defined by
$$Q=\left(\begin{array}{cc}
{2N^2\over k}\sinh(k h(z'))-2h(z') & {2N^2\over k}(\cosh(kh(z'))-1)\\
\noame
\cosh(kh(z'))-1& \sinh(kh(z'))
\end{array}\right)\,.$$
The solution of this system is given by 
\begin{eqnarray}
A(z')=-{h(z')\over 2(1-N^2)}(\partial_{z_1}\bar p(z')-\bar f_1(z'))A_1(z')+ {h(z')\over 2N^2}\bar g_2(z')A_2(z')\,,\nonumber\\
\noame
B(z')=-{h(z')\over 2(1-N^2)}(\partial_{z_1}\bar p(z')-\bar f_1(z'))B_1(z')+ {h(z')\over 2N^2}\bar g_2(z')B_2(z')\,,\nonumber
\end{eqnarray}
where $A_1(z')$, $B_1(z')$ and $A_2(z')$, $B_2(z')$ are solution of
$$Q\left(\begin{array}{c}
A_1\\
\noame
B_1
\end{array}\right)=\left(\begin{array}{c}
h(z')\\
\noame
1
\end{array}\right)\quad\hbox{and}\quad Q\left(\begin{array}{c}
A_2\\
\noame
B_2
\end{array}\right)=\left(\begin{array}{c}
1\\
\noame
0
\end{array}\right).$$
Calculating $A_i$, $B_i$ for $i = 1, 2$, we have 
$$\begin{array}{l} \displaystyle A_1(z')=-{1\over 2},\\
\noame
\displaystyle A_2(z')={\sinh(kh(z'))\over -2h(z')\sinh(kh(z'))+{4N^2\over k}(\cosh(kh(z'))-1)}\,,\\
\noame
 \displaystyle B_1(z')={1\over 2}\coth\left({kh(z')\over 2}\right),\\
 \noame
 \displaystyle B_2(z')={-(\cosh(kh(z')-1)\over -2h(z')\sinh(kh(z'))+{4N^2\over k}(\cosh(kh(z'))-1)}\,,
 \end{array}$$
and then $\bar u_1$ and $\bar w_2$ are obtained by (\ref{sol_1})  as
functions of $\bar p$, $\bar f_1$ and $\bar g_2$, by the following expressions
$$\begin{array}{l}
\bar u_1(z)= \left[{z_3^2\over 2(1-N^2)}+{1\over 4(1-N^2)}\left({2N^2\over k}\sinh(kz_3)-2z_3\right)\right.\\
\noame
\left.-{h(z')\over 2(1-N^2)}{N^2\over k}(\cosh(kh(z'))-1)\coth\left({kh(z')\over 2}\right)\right]\left(\partial_{z_1}\bar p(z')-\bar f_1(z')\right)\\
\noame 
+\Big[-{z_3\over 2N^2}+{h(z')\over 2N^2}\Big(\left({2N^2\over k}\sinh(kz_3)-2z_3\right) A_2\\
\noame
+{2N^2\over k}(\cosh(kz_3)-1)B_2\Big)\Big]\bar g_2(z')\,,
\\
\\
\bar w_2(z)=\Big[{z_3\over 2(1-N^2)}
+{h(z')\over 4(1-N^2)}\Big(\cosh(kz_3)-1\\
\noame
-\coth\left({kh(z')\over 2}\right)\sinh(kz_3)\Big)\Big]\left(\partial_{z_1}\bar p(z')-\bar f_1(z')\right)\\
\noame
 +{h(z')\over 2N^2}\left[\cosh(kz_3)A_2+\sinh(kz_3)B_2\right]\bar g_2(z')\,.
\end{array}$$

As it was pointed at the beginning, expressions for $\bar u_2,\bar w_1$ are  obtained by using the expressions of $\bar u_2$, $\bar w_1$, and so we have
$$\begin{array}{l}
\bar u_2(z)=\left[{z_3^2\over 2(1-N^2)}+{1\over 4(1-N^2)}\left({2N^2\over k}\sinh(kz_3)-2z_3\right)\right.\\
\noame
\left.-{h(z')\over 2(1-N^2)}{N^2\over k}(\cosh(kh(z'))-1)\coth\left({kh(z')\over 2}\right)\right]\left(\partial_{z_2}\bar p(z')-\bar f_2(z')\right)\\
\noame 
-\Big[-{z_3\over 2N^2}+{h(z')\over 2N^2}\left(\left({2N^2\over k}\sinh(kz_3)-2z_3\right)A_2\right.\\
\noame\left.
+{2N^2\over k}(\cosh(kz_3)-1)B_2\right)\Big]\bar g_1(z')\,,
\\
\\
\bar w_1(z)=-\Big[{z_3\over 2(1-N^2)}
+{h(z')\over 4(1-N^2)}\Big(\cosh(kz_3)-1\\
\noame
-\coth\left({kh(z')\over 2}\right)\sinh(kz_3)\Big)\Big]\left(\partial_{z_2}\bar p(z')-\bar f_2(z')\right)\\
\noame
+{h(z')\over 2N^2}\left[\cosh(kz_3)A_2+\sinh(kz_3)B_2\right]\bar g_1(z')\,.
\end{array}$$

We observe $\bar u'$ and $\bar w'$ can be rewritten as follows
\begin{equation}\label{expressions_final_app}\begin{array}{l}
\bar u'(z)= \left[{z_3^2\over 2(1-N^2)}+{1\over 4(1-N^2)}\left({2N^2\over k}\sinh(kz_3)-2z_3\right)\right.\\
\noame 
\left.-{h(z')\over 2(1-N^2)}{N^2\over k}(\cosh(kh(z'))-1)\coth\left({kh(z')\over 2}\right)\right]\left(\nabla_{z'}\bar p(z')-\bar f'(z')\right)\\
\noame 
-\Big[-{z_3\over 2N^2}
+{h(z')\over 2N^2}\left(\left({2N^2\over k}\sinh(kz_3)-2z_3\right)A_2\right.\\
\noame
\left.+{2N^2\over k}(\cosh(kz_3)-1)B_2\right)\Big](\bar g'(z'))^\perp\,,
\end{array}\end{equation}
and
\begin{equation}\label{expressions_final_app123}
\begin{array}{l}
\bar w'(z)=\Big[{z_3\over 2(1-N^2)}+{h(z')\over 4(1-N^2)}\Big(\cosh(kz_3)
 -1\\
\noame
-\coth\left({kh(z')\over 2}\Big)\sinh(kz_3)\right)\Big]\left(\nabla_{z'}\bar p(z')-\bar f'(z')\right)^\perp\\
\noame
 +{h(z')\over 2N^2}\left[\cosh(kz_3)A_2+\sinh(kz_3)B_2\right]\bar g'(z')\,.
\end{array}\end{equation}

Finally, integrating the expressions of $\bar u'$ and $\bar w'$ with respect to the variable $z_3$, it holds that 
\begin{equation}\label{int_u_w_app}
\begin{array}{l}
\displaystyle \int_0^{h(z')}\bar u_j'(z',z_3)\,dz_3=-{h^3(z')\over 1-N^2}\Phi(h(z'),N,R_c)\left(\partial_{z_j}\bar p(z')-\bar f_j(z')\right)\,,\\
\noame
\displaystyle
\int_0^{h(z')}\bar w_j'(z',z_3)\,dz_3=-{1\over 4N^3}\sqrt{{R_c\over 1-N^2}}\Psi(h(y'),N,R_c)\bar g_j(z')\,,
\end{array}\end{equation}
for $j=1,2$, with $\Phi$ and $\Psi$ defined by (\ref{Phi}) and (\ref{Psi}) respectively. Putting this in (\ref{problems_appendix_3}) we get the desired Reynolds equation (\ref{Reynolds_app}).\\


\begin{thebibliography}{10}
\bibitem{Allaire0} Allaire, G.: Homogenization of the Stokes flow in a connected porous medium. Asymp. Anal. 2,  203-222 (1989).

\bibitem{Anguiano_SG} Anguiano, M., Su\'arez-Grau, F.J.: Nonlinear Reynolds equations for non-Newtonian thin-film fluid flows over a rough boundary. IMA J. of Appl. Math. 84,  63-95 (2019).

\bibitem{arbogast} Arbogast, T.,  Douglas J.R., J. , Hornung, U.: Derivation of the double porosity model of single phase flow via homogenization theory. SIAM J. Math. Anal. 21, 823-836 (1990).


\bibitem{Bayada1} Bayada, G., Chambat, M.: The transition between the Stokes equations and the Reynolds equation: a mathematical proof.  Appl.~Math.~Opt. 14, 73-93 (1986).

\bibitem{Bayada_Chambat_2} Bayada, G., Chambat, M.: New models in the theory of the hydrodynamic lubrication of rough surfaces. J. Tribol. 110,  402-407 (1988).

\bibitem{Bayada_Chambat} Bayada, G., Chambat, M.:  Homogenization of the Stokes system in a thin film flow with rapidly varying thickness. RAIRO Mod\'el. Math. Anal. Num\'er. 23, 205-234 (1989). 

\bibitem{BayadaChamGam} Bayada, G., Chambat, M., Gamouana, S.R. : About thin film micropolar asymptotic equations. Quart. Appl. Math. 59, 413-439 (2001).

\bibitem{BCJ} Bayada, G., Ciuperca I., Jai, M.: Homogenized elliptic equations and variational inequalities with oscillating parameters. Application to the study of thin flow behavior with rough surfaces. Nonlinear Anal. Real World Appl. 7, 950-966 (2006). 


\bibitem{BayadaLuc}  Bayada, G., Lukaszewicz, G.:  On micropolar fluids in the theory of lubrication. Rigorous derivation of an analogue of the Reynolds equation. Internat. J. Engrg. Sci. 34, 1477-1490 (1996). 


\bibitem{Benhaboucha} Benhaboucha N., Chambat, M.,  Ciuperca, I.: Asymptotic behaviour of pressure and stresses in a thin film flow with a rough boundary. Quart. Appl. Math. 63, 369-400 (2005).

\bibitem{Benterki} Benterki, D., Benseridi, H.,  Dilmi, M.:  On a non-stationary, non-Newtonian
lubrication problem with Tresca fluid-solid law. J. Inverse Ill-Posed Probl. 27,  719-730 (2019).


\bibitem{Bonnivard}  Bonnivard, M., Pazanin I.,   Su\'arez-Grau, F.J.: Effects of rough boundary and nonzero boundary conditions on the
lubrication process with micropolar fluid. Eur. J. Mech. B Fluids 72,  501-518 (2018).


\bibitem{BP} Boukrouche M., Paoli L.: Asymptotic analysis of a micropolar fluid flow in a thin domain with a free and rough boundary. SIAM J. Math. Anal. 44,  1211-1256 (2012).


\bibitem{BCCM} Bresch, D., Choquet, C., Chupin, L., Colin, T., Gisclon, M.: Roughness-induced effect at main order on the Reynolds approximation. SIAM Multiscale Model. Simul. 8,  997-1017 (2010).

\bibitem{BCiu} Boukrouche, M., Ciuperca, I.: Asymptotic behaviour of solutions of lubrication problem in a thin domain with a rough boundary and Tresca fluid-solid interface law. Quart. Appl. Math. 64,  561-591 (2006).


\bibitem{Ciora} Cioranescu, D., Damlamian, A.,   Griso,  G.: Periodic unfolding and homogenization. C.R. Acad. Sci. Paris Ser. I  335, 99-104 (2002).


\bibitem{Ciora2} Cioranescu, D., Damlamian, A., Griso, G.: The periodic unfolding method in homogenization. SIAM J. Math. Anal.  40,  1585-1620 (2008).

\bibitem{Chupin} Chupin, L.,  Martin, S.: Rigorous derivation of the thin film approximation with roughness-induced correctors. SIAM J. Math. Anal. 44, 2041-3070 (2012).

\bibitem{Dupuy} Dupuy, D., Panasenko, G., Stavre, R.: Asymptotic solution for a micropolar flow in a curvilinear channel.  ZAMM Z. Angew. Math. Mech. 88, 793-807 (2008).

\bibitem{DuvautLions} Duvaut, G., Lions, J.L.: Les inequations en mechanique et en physique [The inequations in mechanics and physics]. Dunod, Paris (1972).

\bibitem{Eringen} Eringen, A.C.: Theory of micropolar fluids. J. Math. Mech. 16, 1-16 (1966).

\bibitem{Fab1} Fabricius, J., Koroleva, Y.O., Tsandzana, A.,  Wall, P.: Asymptotic behaviour of Stokes flow in a thin domain with a moving rough boundary. Proc. R. Soc. A  470: 20130735 (2014).

\bibitem{John} Johnston, G.J., Wayte, R.,   Spikes, H.A.: The measurement and study of very thin lubricant films in concentrated contacts.  Tribol. Trans. 34, 187-194 (1991).

\bibitem{Letoufa} Letoufa, Y., Benseridi, H., Dilmi, M.: Study of Stokes dynamical system in a thin domain with Fourier and Tresca boundary conditions. Asian-European Journal of Mathematics, (2019). https://doi.org/10.1142/S1793557121500078

\bibitem{Luka} Lukaszewicz, G.: Micropolar fluids, theory and applications, Modeling and Simulation
in Science. Birkha$\ddot{\rm u}$ser (1999).

\bibitem{Luo}  Luo, J.B., Huang, P., Wen, S.Z.: Thin film lubrication part I: study on the transition between EHL and thin film lubrication using relative optical interference intensity technique. Wear 194, 107-115 (1996).

\bibitem{Luo2} Luo, J.B., Huang, P., Wen, S.Z., Lawrence,  L.: Characteristics of fluid lubricant films at nano-scale, J. Tribol. 121, 872-878  (1999).

\bibitem{MPM} Marusi\'c-Paloka, E., Pazanin, I.,  Marusi\'c, S.:  An effective model for the lubrication with micropolar fluid. Mech. Res. Comm. 52,  69-73 (2013).


\bibitem{Mikelic2} Mikeli\'c, A.: Remark on the result on homogenization in hydrodynamical lubrication by G. Bayada and M. Chambat.  RAIRO Mod\'el. Math. Anal. Num\'er.  25,  363-370 (1991). 

\bibitem{MT} Mikeli\'c, A., Tapiero, R.: Mathematical derivation of the power law describing polymer flow through a thin slab. Mod\'elisation mathematique et analyse num\'erique 29, 3-21 (1995).

\bibitem{PSG} Pazanin, I.,  Su\'arez-Grau, F.J.: Analysis of the thin film flow in a rough thin domain filled with micropolar fluid. Comput. Math. Appl. 68, 1915-1932  (2014).

\bibitem{Sinha} Singh, C., Sinha, P.: The three-dimensional Reynolds' equation for micropolar fluid lubricated bearings. Wear 76,  199-209 (1982).

\bibitem{grau1} Su\'arez-Grau, F.J.: Asymptotic behavior of a non-Newtonian flow in a thin domain with Navier law on a rough boundary.  Nonlinear Analysis  117, 99-123 (2015).

\bibitem{Tartar} Tartar, L.: Incompressible fluid flow in a porous medium convergence of the homogenization process. In: Appendix to Lecture Notes in Physics, 127. Springer-Velag, Berlin (1980). 
\end{thebibliography}
\end{document}